\documentclass[11pt]{article}
\usepackage{amsmath,amssymb}
\usepackage{epsfig,color}
\usepackage{amsfonts}
\usepackage{booktabs} 
\usepackage{multirow} 
\def\boxit#1{\vbox{\hrule\hbox{\vrule\kern6pt \vbox{\kern6pt#1\kern5pt}
\kern6pt\vrule}\hrule}}

\newcommand{\be}{\utwi{e}}

\newcommand{\by}{\utwi{y}}
\newcommand{\bx}{\utwi{x}}

\newcommand{\bz}{\utwi{z}}
\newcommand{\bu}{\utwi{u}}

\newcommand{\bH}{\utwi{H}}
\newcommand{\bP}{\utwi{P}}

\newcommand{\mN}{\mathbb{N}}
\newcommand{\mR}{\mathbb{R}}

\newcommand{\BX}{\mathbb{X}}

\newcommand{\mF}{{\cal F}}

\newcommand{\mK}{{\cal K}}
\newcommand{\mX}{{\cal X}}

\newcommand{\mL}{{\cal L}}
\newcommand{\mB}{{\cal B}}
\newcommand{\mP}{{\cal P}}

\newcommand{\bmu}{\utwi{\mu}}

\newcommand{\bGamma}{\utwi{\Gamma}}

\newcommand{\tvarsigma}{\tilde{\varsigma}}
\newcommand{\bnu}{\utwi{\nu}}

\newcommand{\bpi}{\utwi{\pi}}

\newcommand{\utwi}[1]{\mbox{\boldmath $#1$}}
\newenvironment{proof}{\trivlist\item[\hskip \labelsep{\sc Proof:}]}
 {\unskip\nobreak\ \lower.3ex\hbox{$\Box$}\endtrivlist}
\newtheorem{theorem}{Theorem}[section]
\newtheorem{lemma}{Lemma}[section]

\def\boxit#1{\vbox{\hrule\hbox{\vrule\kern6pt
          \vbox{\kern6pt#1\kern6pt}\kern6pt\vrule}\hrule}}
\def\JASA{{\it Journal of the American Statistical Association}}

\def\JCGS{{\it Journal of Computational and Graphical Statistics}}

\def\PNAS{{\it Proceedings of the National Academy of Sciences USA}}

\begin{document}
\thispagestyle{empty}

\title{Weak Convergence Rates of Population versus Single-Chain Stochastic Approximation MCMC Algorithms}

\author{Qifan Song,  Mingqi Wu, Faming Liang
 \thanks{correspondence author: Faming Liang.
 Faming Liang is professor, 
 Department of Statistics, Texas A$\&$M University, College Station, TX 77843,
 Email: fliang@stat.tamu.edu;
 Qifan Song is a graduate student,
 Department of Statistics, Texas A$\&$M University, College Station, TX 77843.
 Mingqi Wu is statistical consultant, Shell Global Solutions (US) Inc., 
Shell Technology Center Houston,
3333 Highway 6 South, Houston, TX 77082-3101.
 }}

\maketitle

\begin{abstract}

 In this paper, we establish the theory of weak convergence (toward a normal distribution)
 for both single-chain and population stochastic approximation MCMC algorithms.   
 Based on the theory, we give an explicit ratio of convergence rates for  
 the population SAMCMC algorithm and the single-chain SAMCMC algorithm.  
 Our results provide a theoretic guarantee that the population SAMCMC algorithms 
 are asymptotically more efficient than the single-chain SAMCMC algorithms when the 
 gain factor sequence decreases slower than $O(1/t)$, 
 where $t$ indexes the number of iterations.   
 This is of interest for practical applications.

\vspace{5mm}



\end{abstract}

\newpage

{\centering \section{Introduction}}

 Robbins and Monro (1951) introduced the stochastic approximation algorithm for solving the integration
 equation 
 \begin{equation} \label{Introeq1}
  h(\theta)=\int H(\theta,x) f_{\theta}(x) dx =0,
 \end{equation}
 where $\theta \in \Theta \subset \mR^{d_{\theta}}$ is a parameter vector and $f_{\theta}(x)$, 
 $x \in \mX \subset \mR^{d_x}$, is a density function dependent on $\theta$. The stochastic approximation
 algorithm is a recursive algorithm which proceeds as follows:

\vspace{2mm}

\noindent
 {\it Stochastic Approximation Algorithm}
 \begin{itemize}
 \item[(a)] Draw sample $x_{t+1} \sim f_{\theta_t}(x)$, where $t$ indexes the iteration.
 \item[(b)] Set $\theta_{t+1}=\theta_t+\gamma_{t+1} H(\theta_t,x_{t+1})$, where $\gamma_{t+1}$ is called the 
       gain factor. 
 \end{itemize} 

\vspace{2mm}

 After six decades
 of continual development, this algorithm has developed into an important area in systems control, and 
 has also served as a prototype for development of recursive algorithms for on-line estimation and control
 of stochastic systems. 
 Recently, the stochastic approximation algorithm has been used with Markov chain Monte Carlo (MCMC), 
 which replaces the step (a) by a MCMC sampling step:
 \begin{itemize}
 \item[(a$'$)] Draw a sample $x_{t+1}$ with a Markov transition kernel $P_{\theta_t}(x_t,\cdot)$, which starts with
        $x_t$ and admits $f_{\theta_t}(x)$ as the invariant distribution.
 \end{itemize}
 In statistics, the stochastic approximation MCMC (SAMCMC) algorithm, which is also known as 
 stochastic approximation with Markov state-dependent noise, has been successfully applied to many 
 problems of general interest, such as maximum likelihood estimation
 for incomplete data problems (Younes, 1989; Gu and Kong, 1998), 
 marginal density estimation (Liang, 2007), 
 and adaptive MCMC (Haario {\it et al.}, 2001; Andrieu and Moulines, 2006; 
 Roberts and Rosenthal, 2009; Atchad$\acute{\rm e}$ and Fort, 2009). 
 
 It is clear that efficiency of the SAMCMC algorithm depends 
 crucially on the mixing rate of the Markov transition kernel $P_{\theta_t}$.  
 Motivated by the success of population MCMC algorithms, see e.g., Gilks {\it et al.} (1994),
 Liu {\it et al.} (2000) and Liang and Wong (2000, 2001), which can generally converge 
 faster than single-chain MCMC algorithms, we exploit in this paper the performance of a population 
 SAMCMC algorithm, both theoretically and  numerically. Our results show that the population SAMCMC algorithm 
 can be asymptotically more efficient than the single-chain SAMCMC algorithm. 
 
 Our contribution in this paper is two-fold. First, we establish the asymptotic normality for the  
 SAMCMC estimator, which holds for both the population and single-chain SAMCMC algorithms.
 We note that a similar result has been established in Benveniste {\it et al.} (1990, P.332, Theorem 13),
 but under different conditions for the Markov transition kernel. 
 Our conditions can be easily verified, whereas the conditions given in Benveniste {\it et al.} (1990) are  
 less verifiable.  More importantly, our result is more interpretable than
 that by Benveniste {\it et al.} (1990), 
 and this motivates our design of the population SAMCMC algorithm.
 Second, we propose a general population SAMCMC algorithm, and contrasts its convergence rate with
 that of the single-chain SAMCMC algorithm. Our result provides a theoretical guarantee that the population 
 SAMCMC algorithm is asymptotically more efficient than the single-chain SAMCMC algorithm 
 when the gain factor sequence $\{\gamma_t\}$ decreases slower than $O(1/t)$. 
 The theoretical result has been confirmed with a numerical example. 

 The remainder of this paper is organized as follows. In Section 2, we describe the population SAMCMC algorithm 
 and contrasts its convergence rate with that of the single-chain SAMCMC algorithm.  
 In Section 3, we study the population 
 stochastic approximation Monte Carlo (Pop-SAMC) algorithm, which is proposed based on the 
 SAMC algorithm by Liang {\it et al.}(2007) and is a special case of the population SAMCMC algorithm.  
 In Section 4, we present a numerical example, which compares the performance of SAMC and Pop-SAMC 
 on sampling from a multimodal distribution. 
 In Section 5, we conclude the paper with a brief discussion.

{\centering \section{Convergence Rates of Population versus Single-Chain SAMCMC Algorithms}}

\subsection{Population SAMCMC Algorithm} 

 The population SAMCMC algorithm 
 works with a population of samples at each iteration.  
 Let $\bx_t=(x_t^{(1)}, \ldots, x_t^{(\kappa)})$ denote the population of samples at iteration $t$,
 let $\mX^{\kappa}=\mX \times \cdots \times \mX$ 
 denote the sample space of $\bx_t$, and let $\mX_0^{\kappa}$ denote a subset 
 of $\mX^{\kappa}$ where $\bx_0$ is drawn from. 
 The population SAMCMC algorithm starts with a point $(\theta_0,\bx_0)$ drawn from $\Theta \times \mX_0^{\kappa}$ and 
 then iterates between the following steps:
 
\noindent
 {\it Population Stochastic Approximation MCMC Algorithm} 
 \begin{itemize}
 \item[(a)] Draw samples $x_{t+1}^{(1)}, \ldots, x_{t+1}^{(\kappa)}$ 
       with a Markov transition kernel $\bP_{\theta_t}(\bx_t,\cdot)$, which starts with $\bx_t$ 
       and admits $f_{\theta_t}(\bx)=f_{\theta_t}(x^{(1)}) \ldots f_{\theta_t}(x^{(\kappa)})$ as the invariant distribution.
  \item[(b)] Set $\theta_{t+1}=\theta_t+ \gamma_{t+1} \bH(\theta_t,\bx_{t+1})$, where 
            $\bx_{t+1}=(x_{t+1}^{(1)}, \ldots, x_{t+1}^{(\kappa)})$, and   
            \[  
            \bH(\theta_t, \bx_{t+1}) = \frac{1}{\kappa}\sum_{i=1}^{\kappa} H(\theta_t,x_{t+1}^{(i)}). 
             \]
 \end{itemize}  

  It is easy to see that the population SAMCMC algorithm is actually a SAMCMC algorithm 
  with the mean field function specified by
 \begin{equation} \label{Introeq2}
  h(\theta)= \int \bH(\theta,\bx) \utwi{f}_{\theta}(\bx) d \bx= 
  \int \cdots \int \left [\frac{1}{\kappa} \sum_{i=1}^{\kappa} H(\theta,x^{(i)}) \right] f_{\theta}(x^{(1)}) \ldots 
   f_{\theta}(x^{(\kappa)}) dx^{(1)} \ldots dx^{(\kappa)} =0,
 \end{equation}
 where $\utwi{f}_{\theta}(\bx)=f_{\theta}(x^{(1)}) \ldots f_{\theta}(x^{(\kappa)})$ denotes the joint probability 
 density function of $\bx=(x^{(1)},\ldots, x^{(\kappa)})$.
 
 If $\kappa=1$, the algorithm is reduced to the single-chain SAMCMC algorithm. 
 Compared to the single-chain SAMCMC algorithm, the population SAMCMC algorithm has two 
 advantages. First, it provides 
 a more accurate estimate of $h(\theta)$ at each iteration, and this eventually leads to a faster convergence 
 of the algorithm. 
 Note that $H(\theta_t,\bx_{t+1})$ provides an estimate of $h(\theta_t)$ at iteration $t$. 
 Second, since a population of Markov chains are run in parallel, 
 the population SAMCMC algorithm is able to incorporate some advanced multiple chain operators, such as 
 the crossover operator (Liang and Wong, 2000, 2001), the snooker operator (Gilks {\it et al.}, 1994) and 
 the gradient operator (Liu {\it et al.}, 2000),   
 into simulations. With these operators, the distributed information 
 across the population can then be used in guiding further simulations, and this can 
 accelerate the convergence of the algorithm.  
 However, for illustration purpose,  we consider in this paper primarily the single-chain operator, for which we have 
 \begin{equation} \label{productkernel}
 \bP_{\theta_t}(\bx_t,\bx_{t+1})=\prod_{i=1}^{\kappa} P_{\theta_t}(x_t^{(i)}, x_{t+1}^{(i)} ).
 \end{equation}
 Extension of our convergence result to the general population SAMCMC algorithm 
 which consist of multiple chain operators is straightforward and 
 this will be discussed in Section 3.4.  

\subsection{Main Theoretical Results} 

 For mathematical simplicity, we assume in this paper that $\Theta$ is compact, i.e., the sequence 
 $\{\theta_t\}$ can remain in a compact set. Extension of our results to the case that 
 $\Theta=\mR^{d_{\theta}}$ is trivial with the technique of varying truncations studied in Chen (2002) and   
 Andrieu {\it et al.} (2005), which ensures, almost surely, that  
 the sequence $\{\theta_t\}$ can be included in a compact set.  
 Since Theorems \ref{lem50} and \ref{normalitytheorem} are applicable to 
 both the population and single-chain SAMCMC algorithms, we will let $X_t$ denote the sample(s) 
 drawn at iteration $t$ and let $\BX$ denote the sample space of $X_t$. For the population SAMCMC algorithm,
 we have $\BX=\mX^\kappa$ and $X_t=\bx_t$. For the single-chain SAMC algorithm, we have 
 $\BX=\mX$ and $X_t=x_t$. For any measurable function $f$: $\BX\rightarrow\mathbb{R}^d$, 
 $\bP_{\theta} f(X)=\int_{\BX} \bP_{\theta}(X,y) f(y)d y$. 

 \paragraph{Lyapunov condition on $h(\theta)$.}

 Let $\mL=\{\theta \in \Theta: h(\theta)={\bf 0} \}$.

\begin{itemize}
\item[($A_1$)] The function $h: \Theta \rightarrow \mR^d$ is
 continuous, and there exists a continuously differentiable function $v: \Theta \rightarrow [0,\infty)$
 such that $v_h(\theta)=\nabla^T v(\theta)h(\theta)<0$ for all $\theta \in \mL^c$, 
 $\sup_{\theta \in \mK} v_h(\theta)<0$ for any compact set $\mK \subset \mL^c$,
 and $\nabla v(\theta)$ is Lipschitz continuous.
\end{itemize}

 This condition assumes the existence of a global Lyapunov function $v$ for the mean field $h$.
 If $h$ is a gradient field, i.e., $h=-\nabla J$ for some lower bounded, real-valued and
 differentiable function $J(\theta)$, then $v$ can be set to $J$, provided that $J$ is
 continuously differentiable. This is typical for stochastic optimization problems.

 \paragraph{Stability Condition on $h(\theta)$.}

\begin{itemize}
\item[$(A_2)$] The mean field function $h(\theta)$ is measurable and locally bounded on $\Theta$. 
 There exist a stable matrix $F$ (i.e., all eigenvalues of $F$ are with negative real parts),
 $\rho >0$, and a constant $c$ such that, for any $\theta_* \in \mL$ (defined in $A_1$), 
\[
\| h(\theta)-F(\theta-\theta_*)\| \leq c
\|\theta-\theta_*\|^{2}, \quad \forall \theta \in \{ \theta:
\|\theta-\theta_*\| \leq \rho \}.
\]
\end{itemize}

 This condition constrains the behavior of the mean field function around
 the solution points. 
 If $h(\theta)$ is differentiable, the matrix $F$ can be chosen to be
 the partial derivative of $h(\theta)$, i.e., $\partial h(\theta)/\partial \theta$.
 Otherwise, certain approximation may be needed. 

 \paragraph{Drift condition on the transition kernel $\bP_{\theta}$.}
For a function $g: \BX \rightarrow \mR^{d}$, define the $L_\infty$ norm $\|g\|=\sup_{x \in \BX}\|g(x)\|$,

\begin{itemize}
\item[($A_3$)] For any given $\theta \in \Theta$,
 the transition kernel $\bP_{\theta}$ is irreducible and aperiodic.
 In addition, 
\begin{itemize}
\item[(i)] [Doeblin condition] There exist a constants $\delta>0$, an integer $l>0$ and a probability
 measure $\nu$ such that
\begin{eqnarray}
 & \bullet & \quad \inf_{\theta \in \Theta} \bP^l_{\theta}(X, A) \geq \delta \nu(A), \quad
  \forall X \in \BX, \ \forall A \in \mB_{\BX}, \label{Drieq13}
\end{eqnarray}
 where $\mB_{\BX}$ denotes the Borel set of $\BX$; i.e., the whole set $\BX$ is a {\it small} set for each
 $\bP_\theta$. 

\item[(ii)] There exist a constant $c>0$ such that  for all $X \in \BX$,
\begin{eqnarray}
& \bullet & \quad \sup_{\theta \in \Theta} \|\bH(\theta,\cdot)\| \leq c. \label{Drieq21} \\
& \bullet & \quad \sup_{(\theta,\theta')\in \Theta \times \Theta} \|\theta-\theta'\|^{-1} \|\bH(\theta,\cdot)-\bH(\theta',\cdot)\|
 \leq c. \label{Direq22}
\end{eqnarray}
\item[(iii)]  There exists a constant $c>0$ such that for all $g$ with $\|g\|<\infty$,
\begin{equation}
 \bullet \quad \sup_{(\theta,\theta')\in \Theta\times \Theta}\|\theta-\theta'\|^{-1} 
  \| \bP_{\theta}g-\bP_{\theta'} g\| \leq c \|g \|. \label{Drieq31} \\
\end{equation}
\end{itemize}
\end{itemize}

 The Doeblin condition of Assumption $(A_3)$-(i) is equivalent to assuming that
 the resulting Markov chain has an unique stationary distribution and 
 is uniformly ergodic (Nummelin, 1984, Theorem 6.15). This condition is slightly stronger than the drift condition assumed in
 Andrieu {\it et al.} (2005)  and  Andrieu and Moulines (2006), which implies the $V$-uniform 
 ergodicity for $\bP_\theta$. Assumption $(A_3)$-(ii) gives conditions on $\bH(\theta,X)$, 
 which directly lead to the boundedness of the observation noise.
 It is also worthy to note that the property that $\bP_{\theta}$ satisfies the condition  $(A_3)$-(i)
 and  $(A_3)$-(iii) can be inherited from the corresponding property of the single-chain case. If the conditions hold
 for the single chain kernel $P_\theta$, then the conditions must hold for $\bP_\theta$. One can 
 refer to the arguments used in the proof of Theorem \ref{seesawtheorem} in the supplementary material of this paper.

\paragraph{Conditions on step-sizes.}

\begin{itemize}
\item[$(A_4)$] It consists of two parts:
 \begin{itemize}
 \item[(i)] The sequence $\{\gamma_t\}$, which is defined to be $\gamma(t)$ as a function of $t$ and 
 is exchangeable with $\gamma(t)$ in this paper, 
  is positive and non-increasing and satisfies the following conditions:
 \begin{equation} \label{coneq1}
       \sum_{t=1}^{\infty} \gamma_t=\infty,
  \quad \frac{\gamma_{t+1}-\gamma_t}{\gamma_t}=O(\gamma_{t+1}^{\tau}),
  \quad \sum_{t=1}^\infty \frac{\gamma_t^{(1+\tau')/2}}{\sqrt{t}}<\infty,
\end{equation}
 for some  $\tau \in [1, 2)$ and  $\tau'\in(0,1)$. 
 \item[(ii)] 
 The function $\zeta(t)=\gamma(t)^{-1}$ is differentiable such that its 
 derivative varies regularly with exponent $\tilde{\beta}-1 \geq -1$ 
 (i.e., for any $z>0$, $\zeta'(zt)/\zeta'(t) \to z^{\tilde{\beta}-1}$ 
 as $t \to \infty$), and either of the following two cases holds:
   \begin{itemize}
    \item[(ii.1)] $\gamma(t)$ varies regularly with exponent $(-\beta)$, $\frac{1}{2} < \beta <1$;  

    \item[(ii.2)] For $t \geq 1$, $\gamma(t)=t_0/t$ with $-2 \lambda_F t_0>\max\{1,\tilde{\beta}\}$, where 
     $\lambda_F$ denotes the largest real part of the eigenvalue of the matrix $F$ (defined in condition $A_2$)
     with $\lambda_F<0$.  
 \end{itemize}
\end{itemize}
\end{itemize}

 As shown in Chen (2002, p.134),
 the condition $\sum_{t=1}^\infty \frac{\gamma_t^{(1+\tau')/2}}{\sqrt{t}}<\infty$, together with the monotonicity of 
 $\gamma_t$,  implies   that
 $\gamma_t^{(1+\tau')/2}=o(t^{-1/2})$, and thus 
 \begin{equation} \label{A4derivedeq}
 \sum_{t=1}^\infty \gamma_t^{1+\tau'}=\sum_t (\sqrt{t}\gamma_t^{(1+\tau')/2})(\frac{\gamma_t^{(1+\tau')/2}}{\sqrt{t}})<\infty,
 \end{equation}
 which is often assumed in studying the convergence of stochastic approximations. 
 While condition (\ref{coneq1}) is often assumed in studying the weak convergence of the trajectory averaging estimator of 
 $\theta_t$ (see, e.g., Chen, 2002).
 $(A_4)$-(ii) can be applied to the usual gains $\gamma_t=t_0/t^{\beta}$, $1/2< \beta \leq 1$. 
 Following Pelletier (1998), we deduce that 
 \begin{equation} \label{Peleq1}
 \left( \frac{\gamma_t}{\gamma_{t+1}} \right)^{1/2}= 1+\frac{\beta}{2t} + o(\frac{1}{t}).
 \end{equation}
 In terms of $\gamma_t$, (\ref{Peleq1}) can be rewritten as 
 \begin{equation} \label{Peleq2}
 \left( \frac{\gamma_t}{\gamma_{t+1}} \right)^{1/2}= 1+ \zeta \gamma_t + o(\gamma_t),
 \end{equation} 
 where $\zeta=0$ for the case (ii.1) and $\zeta=\frac{1}{2t_0}$ for $\beta=1$ for the case (ii.2). 
 Clearly, the matrix is $F+\zeta I$ is still stable.

 Theorem \ref{lem50} concerns the convergence of the general
 stochastic approximation MCMC algorithm, whose proof can be found in Appendix A. 
\begin{theorem} \label{lem50}
 Assume that $\Theta$ is compact and the conditions $(A_1)$, $(A_3)$ and $(A_4)$-(i) hold.
 Let the simulation start with a point $(\theta_0,X_0) \in \Theta \times \BX_0$, 
 where $\BX_0 \subset \BX$ such that $\sup_{X \in \BX_0} V(X)<\infty$.  Then, as $t \to \infty$, 
\[
 d(\theta_t,\mL) \to 0, \quad a.s.,
\]
 where $\mL=\{\theta \in \Theta: h(\theta)=0\}$, and $d(u,\bz)=\inf_{z\in \bz} \|u-z\|$. 
\end{theorem}

 To study the convergence rate of $\theta_t$, we rewrite the iterative equation of SAMCMC as 
 \begin{equation} \label{iterequation}
  \theta_{t+1}=\theta_t+\gamma_t [ h(\theta_t)+\xi_{t+1}],
 \end{equation} 
 where $h(\theta_t)=\int_{\BX} \bH(\theta_t, X) f_{\theta_t}(X) dX$,
 and $\xi_{t+1}=\bH(\theta_t, X_{t+1})-h(\theta_t)$ is called the observation noise.
 Lemma \ref{lem51} concerns the decomposition of the observation noise, whose  parts (i) and (iv) are
  partial restatement of Lemma A.5 of Liang (2010). The proof can be found in Appendix B.
 
\begin{lemma} \label{lem51} Assume the conditions of Theorem \ref{lem50} hold. 
  Then there exist $\mR^{d_{\theta}}$-valued random processes $\{e_t\}$,
  $\{\nu_t\}$, and $\{ \varsigma_t \}$ defined on a
  probability space $(\Omega,\mF,\mP)$ such that:
  \begin{itemize}
  \item[(i)] $\xi_t=e_t+\nu_t+\varsigma_t$.
  \item[(ii)] For any constant $\rho>0$ (defined in condition $A_2$), 
         \[
           \begin{split} 
           & E(e_{t+1}|\mF_t)1_{\{\|\theta_t-\theta_*\| \leq \rho \}}=0, \\
           & \sup_{t \geq 0} E( \|e_{t+1}\|^{\alpha}|\mF_t) 1_{\{\|\theta_t-\theta_*\| \leq \rho \}} < \infty, 
           \end{split} 
         \]
    where $\alpha \geq 2$ is a constant and $\mF_t$ is a family of $\sigma$-algebras satisfying
    $\sigma\{\theta_{0}, X_0; \theta_{1}, X_1; \ldots; \theta_t, X_t\} = \mF_t \subseteq \mF_{t+1}$
    for all $t \geq 0$.
  \item[(iii)] Almost surely on $\Lambda(\theta_*)=\{\theta_t \to \theta_*\}$,  as $n\to \infty$, 
             \begin{equation} \label{Gammours}
              \frac{1}{n} \sum_{t=1}^n E(e_{t+1} e_{t+1}'|\mF_t) \to \Gamma, \quad \mbox{a.s.},
              \end{equation}
             where $\Gamma$ is a positive definite matrix. 
  \item[(iv)]$E(\|\nu_{t}\|^2/\gamma_t) 1_{\{\|\theta_t-\theta_*\| \leq \rho \}} \to 0$,  as $t \to \infty$. 
  \item[(v)] $E \|\gamma_t \varsigma_t \| \rightarrow 0$, as $t \rightarrow \infty$.
  \end{itemize}
 \end{lemma}


 This lemma plays a key role in the proof of Theorem \ref{normalitytheorem}, which concerns 
 the asymptotic normality of $\theta_t$. 
 The proof of Theorem \ref{normalitytheorem} can be found in Appendix B. 

 \begin{theorem} \label{normalitytheorem} 
 Assume that $\Theta$ is compact and the conditions $(A_1)$--$(A_4)$ hold.  
 Conditioned on $\Lambda(\theta_*)=\{ \theta_t \to \theta_*\}$, 
\begin{equation} \label{weakeq}
 \frac{\theta_t -\theta_*}{\sqrt{\gamma_t}} \Longrightarrow \mN (0, \Sigma),
 \end{equation}
 with $\Longrightarrow$ denoting the weak convergence, $\mN$ the Gaussian distribution and 
 \begin{equation} \label{Sigmaeq}
 \Sigma=\int_0^{\infty} e^{(F'+\zeta I)t} \Gamma e^{(F+\zeta I)t} dt,
 \end{equation}
 where $F$ is defined in $(A_2)$, $\zeta$ is defined in (\ref{Peleq2}), and 
 $\Gamma$ is defined in Lemma \ref{lem51}. 
\end{theorem} 

\noindent {\bf Remarks}
\begin{itemize} 
 \item[1.] The same result has been established in Benveniste {\it et al.} (1990; Theorem 13, p.332) but under 
  different assumptions for the Markov transition kernel $\bP_{\theta}$. 
  Similar to Andrieu {\it et al.} (2005), we assume a slightly stronger condition $(A_3)$ that $\bP_{\theta}$ satisfies 
  a minorization condition on $\BX$. This condition 
   not only ensures the existence of a stationary distribution of $\bP_{\theta}$, 
  uniform ergodicity,  and the existence and regularity of the solution 
   to the Poisson equation (see e.g., Meyn and Tweedie, 2009), but also implies boundedness of the 
   moment of the sample $X_t$. 
   In Benveniste {\it et al.} (1990), besides some conditions on $\bP_{\theta}$, such as 
   the existence and regularity of the solution to the Poisson equation, 
   the authors impose a moment condition on $X_t$ (Benveniste {\it et al.}, 1990; condition $A_5$, p.220).
   The moment condition is usually very difficult to verify without assumptions on the 
   ergodicity of the Markov chain. 
   Concerning the convergence of the adaptive Markov chain $\{X_t\}$, Andrieu and Moulines (2006) present a central limit theorem
   for the average of $\phi(X_t)$, where $\phi(\cdot)$ is a $V^{r}$-Lipschitz function for some $r \in [0,1/2)$ and 
   $V(\cdot)$ is the drift function.    
   Unlike Andrieu and Moulines (2006),  
   we here present the asymptotic normality for the adaptive stochastic approximation estimator $\theta_t$ itself.

\item[2.] As shown in Benveniste {\it et al.} (1990), 
   $(\theta_t -\theta_*)/\sqrt{\gamma_t}$ converges weakly towards the distribution of a stationary 
   Gaussian diffusion with generator 
 \[ 
  dX_t=(F+\zeta I)X_t+\Gamma^{1/2}dB_t,
  \]  
  where $B_t$ stands for standard Brownian Motion.
  Therefore, the asymptotic covariance matrix $\Sigma$ corresponds to the solution of Lyapunov's equation
 \[
  (F+\zeta I) \Sigma+\Sigma (F'+\zeta I)=- \Gamma.
 \]
 An explicit form of the solution can be found in Ziedan (1972), which is omitted here due to
 its complication.  

 \item[3.]
   From equation (\ref{Gamma}) in the proof of Lemma {\ref{lem51}}, it is not difficult to derive that
  \begin{equation} \label{GammBen}
  \Gamma=\sum_{k=-\infty}^{\infty}\int H(\theta_*, x)[P_{\theta_*}^kH(\theta_*,x)]^Td\pi_{\theta_*}(dx),
   \end{equation}
   where $\pi_{\theta_*}$ denotes the invariant distribution of the transition kernel $P_{\theta_*}$. 
   This is the same expression of $\Gamma$ as given in Benveniste {\it et al.} (1990; equation 4.4.6, p.321).
   Compared to equation (\ref{GammBen}), our expression of $\Gamma$, given in equation (\ref{Gammours}), is more interpretable,
   which corresponds to the asymptotic covariance matrix of $e_t$.
   Given the gain factor sequence  $\{ \gamma_k\}$, the efficiency of a 
   SAMCMC algorithm is determined by $\Gamma$. Based on this observation,
   we show in Theorem \ref{efftheorem} that 
   when $\{\gamma_t\}$ decreases slower than $O(1/t)$, the population SAMCMC algorithm has a smaller 
   asymptotic covariance matrix  
   than the single-chain SAMCMC algorithm and thus is asymptotically more efficient. 

 \item[4.] The condition ``Conditioned on $\Lambda(\theta_*)$'' accommodates the case that there exist 
  multiple solutions for the equation $h(\theta)=0$.  
 \end{itemize}

 Theorem \ref{efftheorem} compares the efficiency of the population SAMCMC and the single-chain 
 SAMCMC, whose proof can be found in Appendix A. 

 \begin{theorem} \label{efftheorem}  
  Suppose that both the population and single-chain SAMCMC algorithms satisfy 
  the conditions  given in Theorem \ref{normalitytheorem}.   
  Let ${\theta}_t^p$ and ${\theta}_t^s$ denote the estimates produced at iteration $t$
  by the population and single-chain SAMCMC algorithms, respectively.
  Given the same gain factor sequence $\{\gamma_t\}$, then  
  $(\theta_t-\theta_*)/\sqrt{\gamma_t}$ and $(\theta_{\kappa t}-\theta_*)/\sqrt{\kappa\gamma_{\kappa t}}$
  have the same asymptotic distribution with the convergence rate ratio
  \begin{equation}\label{effequation}
\frac{\gamma_t}{\kappa\gamma_{\kappa t}}=\kappa^{\beta-1},
  \end{equation}
   where $\kappa$ denotes the population size, and $\beta$ is defined in $(A_4)$. [Note: $1/2 < \beta<1$ for
   the case $A_4$-(ii.1) and $\beta=1$ for the case $A_4$-(ii.2).] 
  \end{theorem}

 \noindent {\bf Remarks}
 \begin{itemize}
  \item[1.] When $\beta=1$ (e.g., $\gamma_t=t_0/t$), the single-chain SAMCMC estimator 
           is as efficient 
          as the population SAMCMC estimator, but this is only true asymptotically.   
          For practical applications, as illustrated by Figure \ref{MSEtimeplot}(a) 
          and Figure \ref{MSEtimeplot50}, the population SAMCMC estimator can still be  
          more efficient than the single-chain SAMCMC estimator due to the population effect: At each iteration, 
          the population SAMCMC provides a more accurate estimate of $h(\theta_t)$ than  
          the single-chain SAMCMC, and this substantially improves the convergence of the algorithm, especially
          at the early stage of the simulation.   

  \item[2.] When $\beta<1$, the population SAMCMC estimator is asymptotically more efficient 
          than the single-chain SAMCMC estimator. This is illustrated by Figure \ref{MSEtimeplot}(b). 
          

  \item[3.] The choice of the population size should be balanced with the choice of $N$, the 
            number of iterations, as the convergence of the algorithm 
            only occurs as $\gamma_t \to 0$. In our experience, 
            $5 \sim 50$ may be a good range for the population size. 
 \end{itemize}

{\centering \section{Population SAMC Algorithm} }

 In this section, we first give a brief review for the SAMC algorithm, and then describe the population 
 SAMC algorithm and its theoretical properties, including convergence and asymptotic normality. 

\subsection{The SAMC Algorithm}

Suppose that we are interested in sampling from a distribution,
\begin{equation} \label{eq1}
  f(x)=c\psi(x), \quad x \in \mX,
\end{equation}
where $\mX$ is the sample space and $c$ is an unknown constant.
Furthermore, we assume that the distribution $f(x)$ is multimodal, which may contain a multitude of 
 modes separated by high energy barriers. 
It is known that the conventional MCMC algorithms, such as the Metropolis-Hastings algorithm 
 (Metropolis {\it et al.}, 1953; Hastings, 1970) and the Gibbs sampler (Geman and Geman, 1984), 
 are prone to get trapped into local modes in simulations from such a kind of distribution.  
 
Designing MCMC algorithms that are immune to the local trap problem has been a long-standing
 topic in Monte Carlo research. A few significant algorithms have been proposed in this direction, 
 including parallel tempering (Geyer, 1991), simulated tempering 
 (Marinari and Parisi, 1992), dynamic weighting (Wong and Liang, 1997), Wang-Landau algorithm 
 (Wang and Landau, 2001), SAMC algorithm (Liang {\it et al.}, 2007), among 
 others.  The SAMC algorithm can be described as follows. 

 Let $E_1,...,E_m$ denote a partition of the sample space $\mX$.  
 For example, the sample space can be partitioned 
 according to the energy function of $f(x)$, i.e., $U(x)=-\log \psi(x)$, into the following 
 subregions: $E_1=\{x: U(x) \leq u_1 \}$, $E_2=\{ x: u_1 < U(x) \leq u_2\}$, 
 $\ldots$, $E_{m-1}=\{x: u_{m-2} < U(x) \leq u_{m-1}\}$ and 
 $E_m=\{x: U(x) \geq u_m\}$, where $u_1 < u_2 < \ldots < u_{m-1}$ are 
 user-specified numbers.  If $\int_{E_i} \psi(x) dx=0$, then $E_i$ is called an empty subregion. 
 Refer to Liang {\it et al.} (2007) 
 for more discussions on sample space partitioning. 
 For the time being, we assume that all the subregions are non-empty; that is,
 $\int_{E_i} \psi(x) dx>0$ for all $i=1,\ldots,m$. 
 Given the partition, SAMC seeks to draw samples from the distribution
\begin{equation} \label{eq2}
  f_w(x)\propto \sum_{i=1}^{m}\frac{\pi_i\psi(x)}{w_i}I(x\in E_i)
\end{equation}
where $w_i=\int_{E_i} \psi(x) dx$, and $\pi_i$'s define the desired sampling frequency
 for each of the subregions and they satisfy the constraints: 
$\pi_i>0$ for all $i$ and $\sum_{i=1}^m \pi_i=1$. If $w_1,...,w_m$ are known,
sampling from $f_w(x)$ will lead to a ``random walk'' in the space of subregions (by regarding each subregion
as a point) with each subregion being sampled with a frequency proportional to $\pi_i$. Thus, the local-trap
problem can be essentially overcome, provided that the sample space is partitioned appropriately.

 Since $w_1, \ldots, w_m$ are generally unknown, SAMC employs the stochastic approximation algorithm 
 to estimate their values. This leads to the following iterative procedure: 

 \noindent {\it The SAMC algorithm}
\begin{itemize}
\item[1.] (Sampling)
 Simulate a sample $x_{t+1}$ by running, for one step, the Metropolis-Hastings algorithm which starts with 
  $x_t$ and admits the stationary distribution: 
\begin{equation} \label{eq5}
  f_{\theta_t}(x)\propto \sum_{i=1}^{m}\frac{\psi(x)}{e^{\theta_{t,i}}}I(x\in E_i),
\end{equation} 
 where $\theta_t=(\theta_{t,1},\ldots, \theta_{t,m})$ and 
 $\theta_{t,i}$ denotes a working (on-line) estimate of $\log(w_i/\pi_i)$ at iteration $t$. 

\item[2.] (Weight updating) Set 
\begin{equation} \label{eq61}
  \theta_{t+1}=\theta_t+\gamma_{t+1} H(\theta_t, x_{t+1}),
\end{equation}
where $H(\theta_t,x_{t+1})=\bz_{t+1}-\bpi$, $\bz_{t+1}=(I(x_{t+1} \in E_1),...,I(x_{t+1} \in E_m))$,  
 $\bpi=(\pi_1, \ldots, \pi_m)$, and $I(\cdot)$ is the indicator function.
\end{itemize}

 A remarkable feature of SAMC is that it possesses the self-adjusting mechanism,  
 which operates based on the past samples. This mechanism
 penalizes the over-visited subregions and rewards the under-visited subregions, and thus enables the system to
 escape from local traps very quickly.
 Mathematically, if a subregion $E_i$ is visited at iteration $t$,  $\theta_{t+1,i}$ will be updated to a larger value,
$\theta_{t+1,i}\leftarrow\theta_{t,i}+\gamma_{t+1}(1-\pi_i)$, such that this subregion has a decreased probability
 to be visited at the next iteration.
On the other hand, for those regions, $E_j \ (j\ne i)$, not visited at iteration $t$, $\theta_{t+1,j}$ will decrease to
a smaller value, $\theta_{t+1,j}\leftarrow\theta_{t,j}-\gamma_{t+1}\pi_j$, such that
 the chance to visit these regions will increase at the next iteration.
 SAMC has been successfully applied to many different problems for which the energy landscape is rugged, 
 such as phylogeny inference (Cheon and Liang, 2009) and Bayesian network learning (Liang and Zhang, 2009).

 \subsection{The Population SAMC Algorithm} \label{ppsamcsection}
 
 The population SAMC (Pop-SAMC) algorithm works as follows. Let $\bx_t=(x_t^{(1)}, \ldots, x_t^{(\kappa)})$ 
 denote the population of samples simulated at iteration $t$. One iteration of the algorithm 
 consists of two steps:

\noindent {\it The Pop-SAMC algorithm:} 
\begin{itemize}
\item[1.] (Population sampling) 
  For $i=1,\ldots, \kappa$, simulate  
  a sample $x_{t+1}^{(i)}$ by running, for one step, the Metropolis-Hasting algorithm  
  which starts with $x_t^{(i)}$ and admits (\ref{eq5}) as the invariant distribution.
  Denote the population of samples by $\bx_{t+1}=(x_{t+1}^{(1)}, \ldots, x_{t+1}^{(\kappa)})$.

\item[2.] (Weight updating) Set 
\begin{equation} \label{eq62}
  \theta_{t+1}=\theta_t+\gamma_{t+1} \bH(\theta_t, \bx_{t+1}), 
\end{equation}
where $\bH(\theta_t,\bx_{t+1})=\sum_{i=1}^{\kappa} H(\theta_t, x_{t+1}^{(i)})/\kappa$, and 
 $H(\theta_t,x_{t+1}^{(i)})$ is as specified in the SAMC algorithm.  
\end{itemize}

 As a special case of the population SAMCMC algorithms, the Pop-SAMC algorithm has a few advantages over the 
 SAMC algorithm. First, since $\bH(\theta,\bx)$ provides a more accurate estimate of $h(\theta)$ 
 than $H(\theta,x)$ at each iteration, Pop-SAMC can converge asymptotically faster than SAMC. 
 This is the so-called  population effect and will be illustrated 
 in Section 4 through a numerical example. 
 Second, population-based proposals, such as the crossover operator, snooker operator  
 and gradient operator, can be included in the algorithm to improve efficiency of 
 the sampling step and thus the convergence of the algorithm.  The only requirement for these operators 
 is that they admit the joint density $f_{\theta_t}(x^{(1)}) \ldots f_{\theta_t}(x^{(\kappa)})$ as the 
 invariant distribution. The weak convergence of the resulting algorithm is discussed at the end of this paper. 
 Third, a smoothing operator can be further introduced to $\bH(\theta,\bx)$ to 
 improve its accuracy as an estimator of $h(\theta)$. Liang (2009) showed through numerical examples that
 the smoothing operator can improve the convergence of SAMC, if multiple MH updates were 
 allowed at each iteration of SAMC.

\subsection{Theoretical Results}

 Regarding the convergence of $\theta_t$, we note that for empty subregions, the corresponding 
 components of $\theta_t$ will trivially converge to $-\infty$ when the number of iterations goes to infinity. 
 Therefore, without loss of generality, we show in the supplementary material only the convergence of the algorithm  
 for the case that all subregions are non-empty; that is, $\int_{E_i} \psi(x) dx>0$ for all $i=1,\ldots,m$.
 Extending the proof to the general case is trivial, since
 replacing (\ref{eq62}) by $(\ref{weiupeq2})$ (given below) will not change the process
 of Pop-SAMC simulation:
\begin{equation} \label{weiupeq2}
\theta_{t+1}'=\theta_t+\gamma_{t+1}(\bH(\theta_t,\bx_{t+1})-\bnu),
\end{equation}
 where $\bnu=(\nu,\ldots,\nu)$ is an $m$-vector of $\nu$, and  
 $\nu=\sum_{j \in \{i: E_i=\emptyset\}} \pi_j/(m-m_0)$ and $m_0$ is the number of empty subregions.

 In our proof, we assume that $\Theta$ is a compact set.
 As aforementioned for the general SAMCMC algorithms,
 this assumption is made only for the reason
 of mathematical simplicity. Extension of our results to the case that
 $\Theta=\mathbb{R}^m$ is trivial with the technique of varying
 truncations (Chen, 2002; Andrieu {\it et al.}, 2005; Liang, 2010). 
 Interested readers can refer to Liang (2010) for the details, where the convergence 
 of SAMC is studied with $\Theta=\mathbb{R}^m$. 
 In the simulations of this paper, we set $\Theta=[-10^{100}, 10^{100}]^m$, as a practical matter,
 this is equivalent to setting $\Theta=\mathbb{R}^m$. 

Under the above assumptions, we have the following theorem concerning the convergence of
  the Pop-SAMC algorithm, whose proof can be found in the supplementary material.

\begin{theorem}  \label{seesawtheorem}
 Let $P_{\theta_t}(x_t^{(i)},x_{t+1}^{(i)})$, $i=1,\ldots, \kappa$ denote the respective 
 Markov transition kernels used for generating the samples $x_{t+1}^{(1)}, \ldots, x_{t+1}^{(\kappa)}$ at 
 iteration $t$. 
 Let $\{\gamma_t\}$ be a gain factor sequence satisfying ($A_4$). 
 If $\Theta$ is compact, all subregions are nonempty,  and  each of the transition kernels satisfies $(A_3)$-$(i)$,  
  then, as $t \to \infty$,   
 \begin{equation} \label{eq7}
 \theta_{t} \rightarrow \theta_*, \quad \mbox{a.s.},
 \end{equation} 
 where $\theta_*=(\theta_*^{(1)}, \ldots, \theta_*^{(m)})$ is given by 
 \begin{equation} \label{eq700} 
  \theta_*^{(i)}= C + \log\left(\int_{E_i} \psi(x) dx\right)-\log(\pi_i), \quad i=1,\ldots, m,
 \end{equation} 
  with $C$ being a constant. 
 \end{theorem}
 The constant $C$ can be determined by imposing a constraint, e.g., $\sum_{i=1}^m e^{\theta_{ti}}$
 is equal to a known number.

\noindent {\bf Remark} As aforementioned, if some regions are empty, the corresponding components of $\theta_*$ will 
 converge to $-\infty$ as $n \to \infty$. In this case, as shown in the supplementary material, we have 
\begin{eqnarray} \label{eq700Sep28}
 \theta_*^{(i)}=
  \begin{cases}
    C + \log\left(\int_{E_i} \psi(x) dx\right)-\log(\pi_i+\nu), & \mbox{if } E_i\ne \emptyset, \\
    -\infty, & \mbox{if } E_i=\emptyset.
  \end{cases}
\end{eqnarray}
 where $C$ is a constant, $\nu=\sum_{j\in\{i:E_i=\emptyset\}} \pi_j/(m-m_0)$,
 and $m_0$ the number of empty subregions. 
 
 The Doeblin condition implies the existence of 
 the stationary distribution $f_{\theta_t}(x)$ for each $\theta_t \in \Theta$, and $P_\theta$ is uniformly ergodic.
 To have this condition satisfied, we assume that $\mX$ is compact and  
 $f(x)$ is bounded away from 0 and $\infty$ on $\mX$. 
 This assumption is true for many Bayesian model selection problems, e.g., change-point
 identification and regression variable selection problems.
 For these problems, after integrating out model parameters from their posterior, 
 the sample space is reduced to a finite set of models.  
 For continuous systems, one may restrict $\mX$ to the region
 $\{x:\psi(x)\geq \psi_{min}\}$, where $\psi_{min}$ is sufficiently small such that the region
 $\{x:\psi(x)<\psi_{min}\}$ is not of interest.  
 For the proposal distribution used in the paper, we assume that it satisfies the local 
 positive condition; that is,  
  there exists two quantities $\epsilon_1>0$ and $\epsilon_2>0$ such that
  $q(x,y)\geq \epsilon_2$ if $|x-y|\leq \epsilon_1$, where $q(x,y)$ denotes the proposal mass/density function. 
  In the supplementary material, we show that the transition kernel induced by local positive proposal satisfies the Doeblin condition.
  The local positive condition is quite standard and has been widely used in the study of MCMC convergence, 
  see, e.g., Roberts and Tweedie (1996). 

 Theorem \ref{SAMCconrateI} concerns the asymptotic normality of $\theta_t$, whose 
 proof can be found in the supplementary material. 

 \begin{theorem} \label{SAMCconrateI} 
 Assume the conditions of Theorem  \ref{seesawtheorem} hold.   
 Conditioned on $\Lambda(\theta_*)=\{\theta_t \to \theta_*\}$, 
\begin{equation} \label{weakeq2}
 \frac{\theta_t -\theta_*}{\sqrt{\gamma_t}} \Longrightarrow \mN (0, \Sigma),
 \end{equation}
 where $\theta_*$ is as defined in (\ref{eq700}), and 
 \[
 \Sigma=\int_0^{\infty} e^{(F'+\zeta I)t} \Gamma e^{(F+\zeta I)t} dt,
 \]
 with $F$ being defined in $(A_2)$, $\zeta$ defined in (\ref{Peleq2}), and
 $\Gamma$ defined in Lemma \ref{lem51}.
\end{theorem}

 Finally, we note that Theorem \ref{efftheorem} is also valid for the SAMC and Pop-SAMC algorithms.  
 Here we would like to emphasize that even when the gain factor sequence 
 is chosen as $\gamma_t=O(1/t)$, Pop-SAMC still has some numerical advantages
 over SAMC in convergence  due to the population effect.
 This will be illustrated by Figure \ref{MSEtimeplot50}.
  
\subsection{Minorization Properties of the Crossover Operator}

 The Pop-SAMC algorithm works on a population of Markov chains.
  Its population setting provides a 
  basis for including more global, advanced MCMC operators, such 
  as the crossover operator of the genetic algorithm, into simulations.  
  Without loss of generality, we assume that the crossover operator works only on the first 
  and second chains of the population. 
 The resulting transition kernel can be written as 
  \begin{equation}\label{cross}
     \bP_{\theta_t}(\bx_t,\bx_{t+1})=P_{\theta_t\times \theta_t}\{(x_t^{(1)},x_t^{(2)}),(x_{t+1}^{(1)},x_{t+1}^{(2)})\}
  \prod_{i=3}^\kappa P_{\theta_t}(x_t^{(i)},x_{t+1}^{(i)}),
  \end{equation}
  which is a product of $\kappa-1$ independent transition kernels, where 
  $P_{\theta_t\times \theta_t}\{(x_t^{(1)},x_t^{(2)}),(x_{t+1}^{(1)},x_{t+1}^{(2)})\}=(1-r_{co})\prod_{i=1}^2 P_{\theta_t}(x_t^{(i)},x_{t+1}^{(i)})
  +r_{co}P_{\theta_t,co}\{(x_t^{(1)},x_t^{(2)}),(x_{t+1}^{(1)},x_{t+1}^{(2)})\}$, and 
  $r_{co}$ is the probability to apply crossover kernel $P_{\theta_t,co}$. 
  Following the proof in the supplementary material, 
  $\prod_{i=1}^2 P_{\theta_t}(x_t^{(i)},x_{t+1}^{(i)})$ is locally positive, 
  which implies that $P_{\theta_t\times \theta_t}$ is locally positive as well, if $r_{co}<1$.
  As long as $\mX$ is compact, and $f(x)$ is bounded away from 0 and $\infty$, $(A_3)$-(i) is satisfied by $P_{\theta_t,co}$.
  The condition $A_3$-(ii) is satisfied because it is independent of the kernel used.   
  The condition $(A_3)$-(iii) can be verified as follows: 
 
  Let $s_{\theta}(\bx,\by)=q(\bx,\by) \min\{1, r(\theta,\bx,\by) \}$, where $\bx=(x^{(1)},x^{(2)})$ and $\by=(y^{(1)},y^{(2)})$, and 
 \[
 r(\theta,\bx,\by)=\frac{f_{\theta}(y^{(1)}) f_{\theta}(y^{(2)})}{f_{\theta}(x^{(1)}) f_{\theta}(x^{(2)})} \frac{q(\by,\bx)}{q(\bx,\by)},
 \]
 is the MH ratio for the crossover operator. It is easy to see that 
 \begin{align*}\left| \frac{\partial s_{\theta}(\bx,\by)}{\partial \theta_i} \right| &= q(\bx,\by) I(r(\theta,\bx,\by) <1) 
  |I(x^{(1)}\in E_i)+I(x^{(2)}\in E_i)-I(y^{(1)}\in E_i)-I(y^{(2)}\in E_i)| r(\theta,\bx,\by)\\
 &\leq2q(\bx,\by).
 \end{align*}
 The mean-value theorem implies that there exists a constant $c$ such that 
 \[
\|s_{\theta}(\bx,\by)-s_{\theta'}(\bx,\by) \| \leq c q(\bx,\by) \| \theta-\theta'\|.
 \]
 Following the same argument as in Liang {\it et al.} (2007),  $(A_3)$-(iii) is satisfied by $P_{\theta_t,co}$. 
 This concludes that each kernel in the right of (\ref{cross}) satisfies the drift condition $(A_3)$.
 Therefore, the product kernel $\bP_{\theta_t}(\bx_t,\bx_{t+1})$ satisfies the drift condition. 
 Then the convergence and asymptotic normality
 of $\theta_t$ (Theorem 3.1 and Theorem 3.2) still hold for this general Pop-SAMC algorithm with crossover operators. 
 We conjecture that the incorporation of crossover operators will bring Pop-SAMC more efficiency. 
 How these advanced operators improve the performance of Pop-SAMC will be explored elsewhere.

{\centering \section{An Illustrative Example} }

To illustrate the performance of Pop-SAMC, 
we study a multimodal example taken from Liang and Wong (2001). The density function over a bivariate
 $\bx$ is given by 
\begin{equation} \label{eq11}
  p(\bx)=\frac{1}{2\pi\sigma^2} \sum_{i=1}^{20}\alpha_i \exp\Big\{-\frac{1}{2\sigma^2}(\bx-\bmu_i)'(\bx-\bmu_i)\Big\},
\end{equation}
where each component has an equal variance $\sigma^2=0.01$ and an equal weight $\alpha_1=...=\alpha_{20}=0.05$, and 
 the mean vectors $\bmu_1, \ldots, \bmu_{20}$ are given in Liang and Wong (2001).  
Since some components of the mixture distribution are far from others, e.g., the distance between the 
 lower right component and its nearest neighboring component is 31.4 times the standard deviation,
 sampling from this distribution puts a great challenge on the existing MCMC algorithms.  


 We set the sample space $\mX=[-10^{100},10^{100}]^2$, and then partitioned it according to the energy 
 function $U(x)=-\log\{p(x)\}$  with an equal energy bandwidth $\Delta u=0.5$ into the 
 following subregions: $E_1=\{x:U(x) \leq 0\}, E_2=\{x:0 < U(x) \leq 0.5\}, ..., E_{20}=\{x:U(x)>9.0\}$.
 Pop-SAMC was first tested on this example with two gain factor sequences, 
 $\gamma_t=100/\max(100,t)$ and $\gamma_t=100/\max(100,t^{0.6})$. 
 In simulations, we set the population size $\kappa=10$, the number of iterations $N=10^6$, and 
 the desired sampling distribution to be uniform, i.e., $\pi_1=\cdots=\pi_{20}=1/20$.  
 The Gaussian random walk proposal distribution was used in the MH sampling step with a 
  covariance matrix of $4I_2$, where $I_2$ is the $2\times 2$ identity matrix. 
 To have a fair comparison with SAMC, we initialize the population in a small region $[0,1]\times [0,1]$,
 which is far from the separated components. 
 Tables \ref{TabGammaI} and \ref{TabGammaII} show the resulting estimates of $P(E_i)$ (i.e. $w_i=\int_{E_i} p(x) dx$), 
  for $i=2,\ldots, 11$, based on 100 independent runs. The computation was done on 
   on a Intel Core 2 Duo 3.0 GHz computer.
 As shown by the true values of $P(E_i)$'s, which are calculated with a total of $2\times 10^9$ samples 
 drawn equally from each of the 20 components of $p(x)$, the subregions $E_2, \ldots, E_{11}$ have 
 covered more than 99\% of the total mass of the distribution.   
 For comparison, SAMC was also applied to this example, but with $N=10^7$ iterations 
 and four gain factor sequences: 
 $\gamma_t=100/\max(100,t)$, $\gamma_t=1000/\max(1000,t)$,
  $\gamma_t=100/\max(100,t^{0.6})$, and $\gamma_t=1000/\max(1000,t^{0.6})$.  
 These settings ensure that each run of Pop-SAMC and SAMC consists of the same number of energy evaluations
 and thus costs about the same CPU time.  
 The resulting estimates of $P(E_i)$'s are summarized in Tables \ref{TabGammaI} and \ref{TabGammaII}.

\begin{table}
\begin{center}
\caption{\label{TabGammaI}
Comparison of efficiency of Pop-SAMC and SAMC for the multimodal example with 
$\gamma_t=t_0/\max\{t_0,t\}$.
The number in the parentheses shows the standard error of the estimate of $P(E_i)$.
 }
\begin{tabular}{ccccccc} \\ \hline
            &  & Pop-SAMC & & SAMC  & &  SAMC \\ \cline{3-3} \cline{5-5} \cline{7-7}
Setting     & True         & $(t_0,\tau,N)=(100,10,10^6)$ &&  $(t_0,N)=(100,10^7)$  & & $(t_0,N)=(1000,10^7)$ \\ \hline
$P(E_2)$    & 0.2387    & 0.2383(0.0003) & & 0.2390(0.0003) & &  0.2382(0.0008)\\
$P(E_3)$    & 0.3027    & 0.3027(0.0003) & & 0.3024(0.0003) & &  0.3030(0.0008)\\
$P(E_4)$    & 0.1856    & 0.1859(0.0002) & & 0.1859(0.0002) & &  0.1852(0.0006)\\
$P(E_5)$    & 0.1124    & 0.1124(0.0001) & & 0.1121(0.0001) & &  0.1126(0.0004)\\
$P(E_6)$    & 0.0663    & 0.0663(0.0001) & & 0.0662(0.0001) & &  0.0666(0.0003) \\
$P(E_7)$    & 0.0384    & 0.0384(0)      & & 0.0384(0)      & &  0.0383(0.0001) \\
$P(E_8)$    & 0.0226    & 0.0226(0)      & & 0.0225(0)      & &  0.0227(0.0001) \\
$P(E_9)$    & 0.0134    & 0.0134(0)      & & 0.0134(0)      & &  0.0135(0.0001)\\
$P(E_{10})$ & 0.0080    & 0.0080(0)      & & 0.0080(0)      & &  0.0079(0)  \\
$P(E_{11})$ & 0.0048    & 0.0048(0)      & & 0.0048(0)      & &  0.0048(0) \\ \hline
 CPU (s)    & ---       & 18           & &  21          & &  21       \\ \hline
\end{tabular}
\end{center}
\end{table}

 Our numerical results agree extremely well with Theorem \ref{efftheorem}. 
 It follows from the delta method (see, e.g., Casella and Berger, 2002) that the mean square errors (MSEs) 
 of the estimates of $P(E_i)$ should follow the same limiting rule (\ref{effequation}) as $\theta_t$ does. 
  For this example, when the same gain factor sequence $\gamma_t=100/\max(100,t)$ is used, 
  SAMC is as efficient as Pop-SAMC when the number of iterations is large; 
  the two estimators share the same standard errors as reported in Table \ref{TabGammaI}. 
  When the gain factor sequence $\gamma_t=1000/\max(1000,t)$ is used for SAMC, 
  the runs of SAMC and Pop-SAMC end 
  with the same gain factor values. In this case, as expected, the SAMC estimator has  
  larger standard errors than the Pop-SAMC estimator; 
  the relative efficiency of these two estimators is about $3.0^2$ 
  ($3 \approx (0.0008+\cdots+0.0003)/(0.0003+\cdots+0.0001)$),  which is close to the theoretical value 10.  
  When the gain factor sequence $\gamma_t=100/\max(100,t^{0.6})$ is used,
  Pop-SAMC is more efficient than SAMC. Table \ref{TabGammaII} shows that 
  the relative efficiency of the Pop-SAMC estimator versus the SAMC estimator is about 2.56 ($=1.6^2$ 
   and $1.6\approx (0.0065+\cdots+0.0002)/(0.0042+\cdots+0.0001)$), 
   which agrees well with the theoretical value 2.51 ($=10^{0.4}$).

\begin{table}
\begin{center}
\caption{\label{TabGammaII}
Comparison of efficiency of Pop-SAMC and SAMC for the multimodal example with $\gamma_t=t_0/\max\{t_0,t^{0.6}\}$.
 The number in the parentheses shows the standard error of the estimate of $P(E_i)$.
 }
\begin{tabular}{ccccccc} \\ \hline
            &  & Pop-SAMC & & SAMC  & &  SAMC \\ \cline{3-3} \cline{5-5} \cline{7-7}
Setting     & True         & $(t_0,\tau,N)=(100,10,10^6)$ &&  $(t_0,N)=(100,10^7)$  & & $(t_0,N)=(1000,10^7)$ \\ \hline
$P(E_2)$    & 0.2387    & 0.2236(0.0042) &&   0.2244(0.0065) &&   0.1534(0.0184) \\
$P(E_3)$    & 0.3027    & 0.3045(0.0041) &&   0.3123(0.0076) &&   0.3329(0.0268) \\
$P(E_4)$    & 0.1856    & 0.1909(0.0035) &&   0.1850(0.0054) &&   0.1815(0.0207) \\
$P(E_5)$    & 0.1124    & 0.1156(0.0024) &&   0.1167(0.0039) &&   0.1205(0.0144) \\
$P(E_6)$    & 0.0663    & 0.0706(0.0015) &&   0.0648(0.0020) &&   0.0715(0.0096) \\
$P(E_7)$    & 0.0384    & 0.0387(0.0008) &&   0.0390(0.0013) &&   0.0542(0.0071) \\
$P(E_8)$    & 0.0226    & 0.0227(0.0005) &&   0.0243(0.0010) &&   0.0303(0.0058) \\
$P(E_9)$    & 0.0134    & 0.0133(0.0003) &&   0.0137(0.0005) &&   0.0218(0.0069) \\
$P(E_{10})$ & 0.0080    & 0.0082(0.0002) &&   0.0079(0.0003) &&   0.0184(0.0050) \\
$P(E_{11})$ & 0.0048    & 0.0047(0.0001) &&   0.0047(0.0002) &&   0.0047(0.0009) \\ \hline
 CPU(s)    & ---       & 18           &&   23           &&   22       \\ \hline
\end{tabular}
\end{center}
\end{table}

\begin{figure}[htbp]
\begin{center}
\begin{tabular}{cc}
(a) & (b) \\
\epsfig{figure=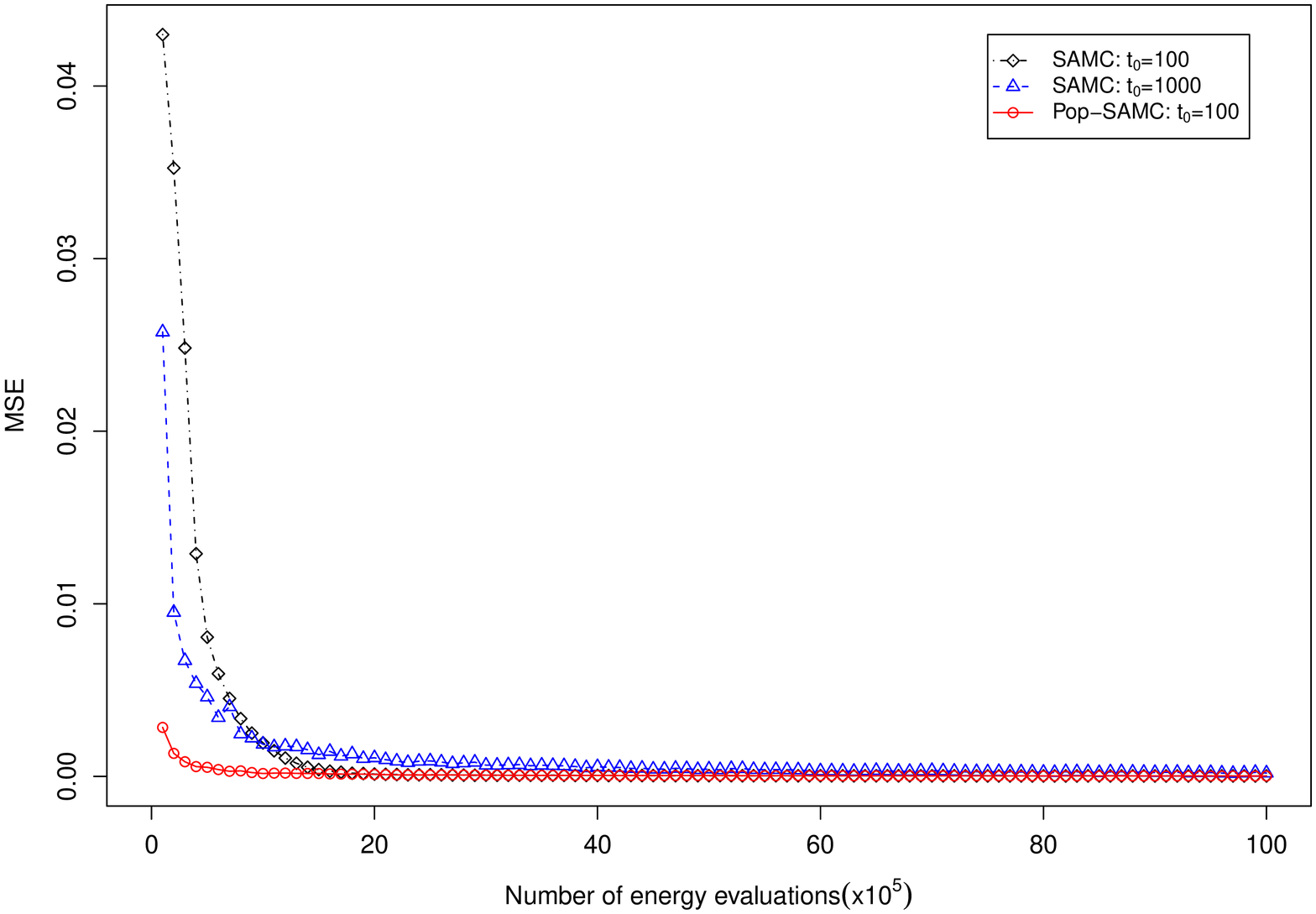, height=3.0in,width=2.5in, angle=270}
 & 
\epsfig{figure=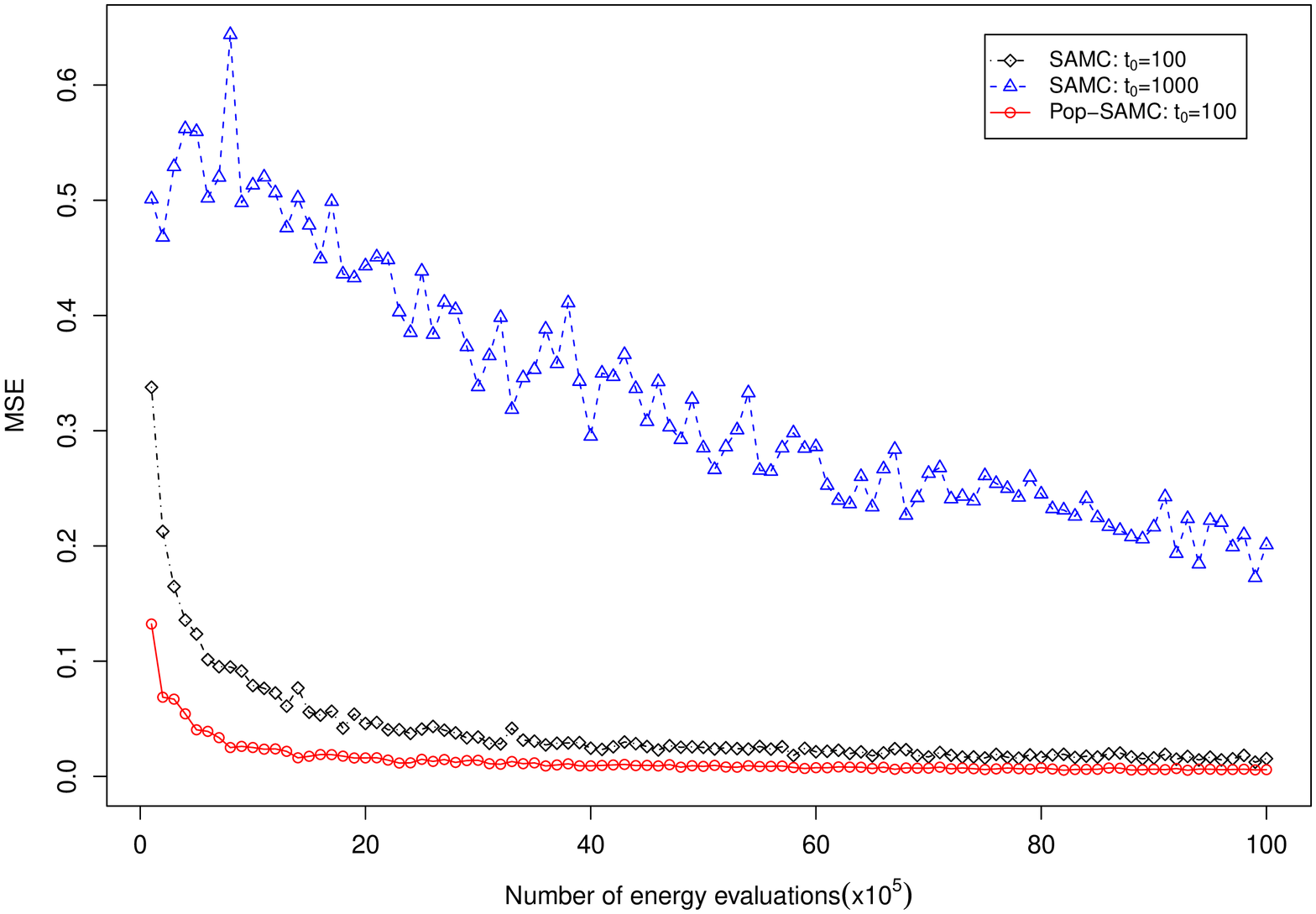, height=3.0in,width=2.5in, angle=270}
\end{tabular}
\caption{ \label{MSEtimeplot}
Mean square errors (MSEs) produced by Pop-SAMC and SAMC at different iterations. The left plot is
produced with  $\gamma_t=t_0/\max(t_0,t)$, and the right plot is
produced with  $\gamma_t=t_0/\max(t_0,t^{0.6})$. }
\end{center}
\end{figure}

 We note that the results reported in Tables \ref{TabGammaI} and \ref{TabGammaII} are only 
 for the scenario that the number of iterations is large. For a thorough comparison, 
 we evaluated the MSEs of the Pop-SAMC and SAMC estimators at 100 equally spaced time points, 
 with iterations  $10^4 \sim 10^6$ for Pop-SAMC and $10^5 \sim 10^7$ for SAMC.  
 The results are shown in Figure \ref{MSEtimeplot}. 
 The plots indicate that Pop-SAMC can converge much faster than SAMC, even when the gain 
 factor sequence $\gamma_t=t_0/\max(t_0,t)$ is used. As discussed previously, this is due to the 
 population effect: Pop-SAMC provides a more accurate estimator of 
 $h(\theta_t)$ at each iteration, and this improves its convergence, 
 especially at the early stage of the simulation. 
 
\begin{figure}[htbp]
\begin{center}
\begin{tabular}{c}
\epsfig{figure=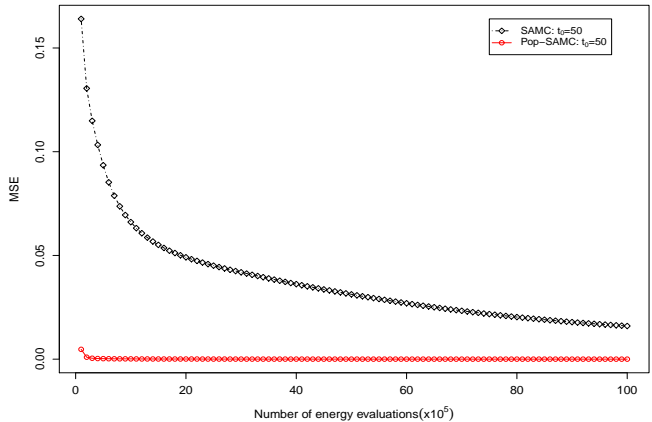, height=3.5in,width=2.5in, angle=270}
\end{tabular}
\caption{ \label{MSEtimeplot50}
Mean square errors (MSEs) produced by Pop-SAMC and SAMC at different iterations with the gain factor 
sequence $\gamma_t=50/\max(50,t)$. }
\end{center}
\end{figure}

 To further explore the population effect of Pop-SAMC, both Pop-SAMC and SAMC were re-run 
 100 times with a smaller gain factor sequence $\gamma_t=50/\max(50,t)$.  
 Figure \ref{MSEtimeplot50} shows that under this setting, 
 SAMC converges very slowly, while Pop-SAMC still converges very fast. 
 This experiment shows that Pop-SAMC is more robust to the choice of gain factor sequence, 
 and it can work with a smaller gain factor sequence than can SAMC. 

{\centering \section{Conclusion} }

 In this paper, we have proposed a population SAMCMC algorithm and contrasted its 
 convergence rate with that of the single-chain SAMCMC algorithm. As the main theoretical result, 
 we establish the limiting ratio between the $L_2$ rates of convergence of 
 the two types of SAMCMC algorithms. Our result provides a theoretical guarantee that the population SAMCMC 
 algorithm is asymptotically more efficient than the single-chain SAMC algorithm 
 when the gain factor sequence $\{\gamma_t\}$ decreases slower than $O(1/t)$.  
 This theoretical result has been confirmed with a numerical example. 

 In this paper, we have also proved the asymptotic normality of SAMCMC estimators under 
 mild conditions. As mentioned previously, the major difference between 
 this work and Benveniste {\it et al.} (1990) are the assumptions on 
 Markov transition kernels. Our assumptions are easier 
 to verify than those by Benveniste {\it et al.} (1990).  
 We note that the work by Chen (2002) and Pelletier (1998) can potentially be
 extended to SAMCMC algorithms. The major differences between their work and ours  
 are the assumptions on observation noise. In Chen (2002) (Theorem 3.3.2, p.128) 
 and Pelletier (1998), it is assumed that the observation noise 
 can be decomposed in the form
\[
 \epsilon_t=e_t+\nu_t,
\]
 where $\{e_t\}$ forms a martingale difference sequence and  $\{\nu_t\}$ is a higher order 
 term of $O(\sqrt{\gamma_t})$. However, as shown in 
 Lemma \ref{lem51}, the SAMCMC algorithms do not satisfy this assumption. 


   

\vspace{1cm}

{\centering \section*{Acknowledgment} }

 Liang's research was supported in part by the grant (DMS-1106494 and DMS-1317131) made by
 the National Science
 Foundation and the award (KUS-C1-016-04) made by King Abdullah University of
 Science and Technology (KAUST). 
 The authors thank the editor,
 the associate editor, and the referee for their constructive comments which have led to
 significant improvement of this paper.

\appendix

{\centering \section*{Appendix}} 

\section{Proof of Theorem \ref{lem50} } 

 To prove Theorem \ref{lem50}, we first introduce the following lemmas.
 Lemma \ref{lem1} is a combined restatement of Theorem 2 of Andrieu and Moulines (2006), 
 Proposition 6.1 of Andrieu {\it et al.} (2005), and Lemma 5 of Andrieu and Moulines (2006).

 \begin{lemma} \label{lem1}
 Assume that $\Theta$ is compact and the condition $(A_3)$ holds. Then the following results hold:
 \begin{itemize}
 \item[$(B_1)$] For any $\theta \in \Theta$, the Markov kernel $P_{\theta}$ has a single
 stationary distribution $\pi_{\theta}$. 
 In addition, $H: \Theta \times \BX \to \Theta$ is
 measurable and for all $\theta \in \Theta$, $\int_{\BX} \|H(\theta,x)\| \pi_{\theta}(x) dx <\infty$.

 \item[$(B_2)$] For any $\theta \in \Theta$, the Poisson equation $u_{\theta}(X)-P_{\theta}u_{\theta}(X)
 =H(\theta,X)-h(\theta)$ has a  solution $u_{\theta}(X)$, where $P_{\theta} u_{\theta}(X)=
 \int_{\BX} u_{\theta}(y) P_{\theta}(X,y) d y$.
 For any $\eta\in(0,1)$,
 the following conditions hold: 
 \begin{equation} \label{boundeq1}
 \begin{split}
 (i)& \quad \sup_{\theta \in \Theta}\left( \|u_{\theta}(\cdot)\|+\|P_{\theta} u_{\theta}(\cdot)\| \right) < \infty, \\
 (ii)& \quad  \sup_{(\theta,\theta')\in \Theta \times \Theta} \|\theta-\theta'\|^{-\eta} 
 \left\{ \|u_{\theta}(\cdot)-u_{\theta'}(\cdot)\|+\|P_{\theta} u_{\theta}(\cdot)-P_{\theta'} u_{\theta'}(\cdot)\| \right\} <\infty. \\
 \end{split}
\end{equation}

\item[$(B_3)$] For any $\eta \in (0, 1)$, 
    \[
    \sup_{ (\theta,\theta') \in \Theta \times \Theta } \| \theta-\theta'\|^{-\eta} \|h(\theta)-h(\theta')\| < \infty. 
   \]
\end{itemize}
\end{lemma}

 
 Tadi$\acute{\mbox{c}}$ (1997) studied the convergence of the stochastic approximation MCMC algorithm
 under different conditions from those given in Andrieu, Moulines and Priouret (2005) and Andrieu and Moulines (2006).
 We combined some results of the three papers and got the following lemma, which corresponds to 
 Theorem 4.1 and Lemma 2.2 of Tadi$\acute{\mbox{c}}$ (1997). 
 
\begin{lemma} \label{lem2} Assume the conditions of Theorem \ref{lem50} hold. Then the following results hold:
\begin{itemize}
\item[$(C_1)$] There exist $\mR^{d_{\theta}}$-valued random processes $\{\epsilon_t\}_{t \geq 0}$,
 $\{\epsilon_t^{'}\}_{t \geq 0}$ and $\{\epsilon_t^{''}\}_{t \geq 0}$ defined on a
 probability space $(\Omega,\mF,\mP)$ such that
 \begin{equation} \label{seleq1}
 \gamma_{t+1} \xi_{t+1}=\epsilon_{t+1}+\epsilon_{t+1}'+\epsilon_{t+1}^{''}-\epsilon_t^{''}, \quad t \geq 0,
 \end{equation}
 where $\xi_{t+1}=H(\theta_t,X_{t+1})-h(\theta_t)$. 
\item[$(C_2)$] The series
 $\sum_{t=0}^{\infty} \| \epsilon_t' \|$,  $\sum_{t=0}^{\infty} \|\epsilon_t^{''}\|^2$ and
 $\sum_{t=0}^{\infty} \|\epsilon_{t+1}\|^2$ all converge a.s. and
\begin{equation} \label{seleq2}
  E(\epsilon_{t+1} |\mF_t)=0, \quad a.s., \quad n \geq 0,
\end{equation}
where $\{\mF_t\}_{t \geq 0}$ is a family of $\sigma$-algebras of $\mF$ satisfying $\sigma\{\theta_0\} \subseteq
 \mF_0$ and $\sigma\{\epsilon_t, \epsilon_t', \epsilon_t^{''} \} \subseteq \mF_t \subseteq \mF_{t+1}$, $t \geq 0$.

\item[$(C_3)$] Let $R_t=R_t'+R_t^{''}$, $t \geq 1$, where $R_t'=\gamma_{t+1} \nabla^T v(\theta_t) \xi_{t+1}$, and
 \[
R_{t+1}^{''}=\int_0^1 \big[\nabla v(\theta_t+s(\theta_{t+1}-\theta_t))-\nabla v(\theta_t)\big]^T
 (\theta_{t+1}-\theta_t) ds.
\]
 Then $\sum_{t=1}^{\infty} \gamma_t \xi_t$ and $\sum_{t=1}^{\infty} R_t$ converge a.s..
\end{itemize}
\end{lemma}
 
\begin{proof} \begin{itemize}
 \item[$(C_1)$]  
 Let $\epsilon_0=\epsilon_0'=0$, and
\[
\begin{split}
\epsilon_{t+1} &=\gamma_{t+1} \big[ u_{\theta_t}(x_{t+1})-P_{\theta_t} u_{\theta_t}(x_t) \big], \\
\epsilon_{t+1}' &=\gamma_{t+1} \big[ P_{\theta_{t+1}}u_{\theta_{t+1}}(x_{t+1})-
 P_{\theta_{t}}u_{\theta_{t}}(x_{t+1})\big] +(\gamma_{t+2}-\gamma_{t+1})P_{\theta_{t+1}}u_{\theta_{t+1}}(x_{t+1}), \\
 \epsilon_t^{''} &=-\gamma_{t+1} P_{\theta_{t}}u_{\theta_{t}}(x_{t}). \\
\end{split}
\]
It is easy to verify that (\ref{seleq1}) is satisfied.

\item[$(C_2)$] Since 
\[
E(u_{\theta_t}(x_{t+1})|\mF_t)= P_{\theta_t}u_{\theta_t}(x_t),
\]
which concludes (\ref{seleq2}).  It follows from $(B_2)$, $(A_3)$ and $(A_4)$ that there exist
 constants $c_3,c_4,c_5,c_6,c_7 \in \mR^+$ such that 
\[
\begin{split}
\|\epsilon_{t+1}\|^2 & \leq 2 c_3 \gamma_{t+1}^2,\quad
\|\epsilon_{t+1}^{''}\|^2 \leq c_4 \gamma_{t+1}^2 , \\
\| \epsilon_{t+1}'\| & \leq c_5 \gamma_{t+1}\|\theta_{t+1}-\theta_t\|^\eta
 + c_6 \gamma_{t+1}^{1+\tau} \leq c_7 \gamma_{t+1}^{1+\eta},\\
\end{split}
\]
for any $\eta\in(0,1)$. 
 Following from  (\ref{A4derivedeq}) and setting $\eta\geq\tau'$ ($\tau'$ is defined in $A_4$), we have 
\[
\sum_{t=0}^{\infty} \|\epsilon_{t+1}\|^2  <\infty, \quad \sum_{t=0}^{\infty}  \| \epsilon_{t+1}'\| <\infty, 
\quad  \sum_{t=0}^{\infty} \|\epsilon_{t+1}^{''}\|^2 <\infty, 
\]
which, by Fubini's theorem, implies that the series
 $\sum_{t=0}^{\infty} \|\epsilon_{t+1}\|^2$,  $\sum_{t=0}^{\infty}  \| \epsilon_{t+1}'\|$, and 
 $\sum_{t=0}^{\infty} \|\epsilon_{t+1}^{''}\|^2$ all converge almost surely to some finite value random variables.  

\item[$(C_3)$] Let $M=\sup_{\theta \in \Theta} \max\{\|h(\theta)\|, \| \nabla v(\theta) \| \}$,
 and $L$ is the Lipschitz constant of $\nabla v(\cdot)$.
 Since $\sigma\{\theta_t\} \subset \mF_t$, it follows from $(C_2)$ that
 $E(\nabla^T v(\theta_t) \epsilon_{t+1} |\mF_t)=0$. In addition, we have
\[
\sum_{t=0}^{\infty} E \big( |\nabla^T v(\theta_t) \epsilon_{t+1}| \big)^2 \leq
 M^2 \sum_{t=0}^{\infty} E \big( \| \epsilon_{t+1} \|^2  \big) <\infty.
\]
It follows from the martingale convergence theorem (Hall and Heyde, 1980; Theorem 2.15)
 that both $\sum_{t=0}^{\infty} \epsilon_{t+1}$ and
 $\sum_{t=0}^{\infty} \nabla^T v(\theta_t) \epsilon_{t+1}$ converge almost surely. Since
\[
\begin{split}
\sum_{t=0}^{\infty} |\nabla^T v(\theta_t) \epsilon_{t+1}' | & \leq M \sum_{t=1}^{\infty} \| \epsilon_t'\|, \\
 \sum_{t=1}^{\infty} \gamma_t^2 \|\xi_t\|^2 & \leq C\left( \sum_{t=1}^{\infty} \|\epsilon_t\|^2
+ \sum_{t=1}^{\infty} \|\epsilon_t'\|^2+\sum_{t=0}^{\infty} \|\epsilon_t^{''}\|^2\right), \\
\end{split}
\]
for some constant $C$. It follows from $(C_2)$ that both
$\sum_{t=0}^{\infty} |\nabla^T v(\theta_t) \epsilon_{t+1}' |$ and
$\sum_{t=1}^{\infty} \gamma_t^2 \|\xi_t\|^2$ converge. In addition,
\[
\begin{split}
 &\|R_{t+1}^{''}\| \leq L \| \theta_{t+1}-\theta_t\|^2 =
  L \| \gamma_{t+1} h(\theta_t)+\gamma_{t+1} \xi_{t+1} \|^2
 \leq  2L \big( M^2 \gamma_{t+1}^{2}+ \gamma_{t+1}^2 \|\xi_{t+1}\|^2 \big),  \\
& \left|\left( \nabla v(\theta_{t+1})-\nabla v(\theta_t) \right)^T \epsilon_{t+1}^{''} \right|
 \leq L \|\theta_{t+1}-\theta_t\| \| \epsilon_{t+1}^{''} \|, \\
\end{split}
\]
for all $t \geq 0$. Consequently,
\[
\begin{split}
\sum_{t=1}^{\infty} |R_t^{''}| & \leq 2 L M^2 \sum_{t=1}^{\infty} \gamma_t^2+2 L \sum_{t=1}^{\infty}
\gamma_t^2 \|\xi_t\|^2 <\infty, \\
\sum_{t=0}^{\infty} \left| \left(\nabla v(\theta_{t+1})-\nabla v(\theta_t)\right)^T \epsilon_{t+1}^{''} \right| &\leq
 \left( 2L^2 M^2 \sum_{t=1}^{\infty} \gamma_t^2 +2 L^2 \sum_{t=1}^{\infty} \gamma_t^2 \|\xi_t\|^2 \right)^{1/2}
 \left( \sum_{t=1}^{\infty} \|\epsilon_t^{''}\|^2 \right)^{1/2} <\infty. \\
\end{split}
\]
Since
\[
\begin{split}
\sum_{t=1}^{n} \gamma_t \xi_t &=\sum_{t=1}^{n} \epsilon_t +\sum_{t=1}^{n} \epsilon_t' +\epsilon_n^{''}
 -\epsilon_0^{''}, \\
\sum_{t=0}^{n} R_{t+1}' &= \sum_{t=0}^n \nabla^T v(\theta_t) \epsilon_{t+1}+
 \sum_{t=0}^{n} \nabla^T v(\theta_t) \epsilon_{t+1}'-\sum_{t=0}^{n}
 \left(\nabla v(\theta_{t+1})-\nabla v(\theta_t) \right)^T \epsilon_{t+1}^{''} \\
   & + \nabla^T v(\theta_{n+1}) \epsilon_{n+1}^{''}-\nabla^T v(\theta_0) \epsilon_0^{''}, \\
\end{split}
\]
and $\epsilon_n''$ convergent to zero by $(C_2)$,  it is obvious that $\sum_{t=1}^{\infty} \gamma_t \xi_t$ and $\sum_{t=1}^{\infty} R_t$ converge almost surely.
\end{itemize}
 The proof for Lemma \ref{lem2} is completed.
\end{proof}

 Based on Lemma \ref{lem2}, Theorem \ref{lem50} can be proved in a similar way to Theorem 2.2 of 
 Tadi$\acute{\mbox{c}}$ (1997). Since Tadi$\acute{\mbox{c}}$ (1997) is not available publicly, 
 we reproduce the proof for Theorem \ref{lem50} in Supplemental Materials.

\section{Proofs of Lemma \ref{lem51}, Theorem \ref{normalitytheorem} and Theorem \ref{efftheorem} }

\subsection{Proof of Lemma \ref{lem51}.} 
 
 Lemma \ref{lem60} is a restatement of Proposition 6.1 of Andrieu {\it et al.} (2005). It has a little overlap with $(B_2)$. 

\begin{lemma} \label{lem60} Assume $A_3$-(i) and $A_3$-(iii) hold. Suppose that the family of functions $\{ g_{\theta}, \theta \in \Theta\}$ 
 satisfies the condition:   For any compact subset $\mK\subset\Theta$,
\begin{equation} \label{fm2012eq1}
 \sup_{\theta\in\mathcal K}\|g_\theta(\cdot)\|<\infty,\qquad
 \sup_{(\theta,\theta')\in\mathcal K\times\mathcal K}|\theta-\theta'|^{-\iota}\|g_{\theta}(\cdot)-g_{\theta'}(\cdot)\|<\infty,
\end{equation} 
 for some  $\iota\in (0,1)$. Let $u_{\theta}(x)$ be the solution to the Poisson equation 
 $u_{\theta}(x)-P_{\theta} u_{\theta}(x)=g_{\theta}(x)-\pi_{\theta}(g_{\theta}(x))$, 
 where $\pi_{\theta}(g_{\theta}(x))=\int_{\BX} g_{\theta}(x) \pi_{\theta}(x) dx$. 
 Then, for any compact set $\mK$ and any $\iota'\in(0,\iota)$, 
\[
 \begin{split} 
 \sup_{\theta\in\mK}& \left(\|u_{\theta}(\cdot)\| +\| P_{\theta} u_{\theta}(\cdot)\| \right) < \infty, \\
 \sup_{(\theta,\theta') \in \mK \times \mK}&\| \theta-\theta'\|^{-\iota'} \left\{\|u_{\theta}(\cdot)-u_{\theta'}(\cdot)\| 
   +\| P_{\theta} u_{\theta}(\cdot) - P_{\theta'} u_{\theta'}(\cdot)\| \right\} <  \infty. \\ 
\end{split} 
 \]
 \end{lemma}

 Lemma \ref{slln} can be viewed as a partial restatement of Proposition 7 of Andrieu and Moulines (2006), but under 
 different conditions. 

\begin{lemma}\label{slln}
 Assume that $\Theta$ is compact and the conditions ($A_3$) and ($A_4$)-(i) hold.  
 Let $\{g_\theta, \theta\in\Theta\}$ be a family of 
 functions satisfying (\ref{fm2012eq1}) 
 with $\iota\in((1+\tau')/2,1)$, where $\tau'$ is defined in condition $A_4$.
  Then
\[
n^{-1}\sum_{k=1}^n\left(g_{\theta_k}(X_k)-\int_{\BX} g_{\theta_k}(x)d\pi_{\theta_k}(x)\right)\rightarrow 0, \qquad a.s.
\] 
 for any starting point $(\theta_0,X_0)$. 
\end{lemma}
\begin{proof}
Without loss of generality, we assume that $g_\theta$ takes values on $\mathbb{R}$.
(If $g_{\theta}$ takes vales on $\mathbb{R}^d$, the proof can be done elementwisely.) 
Let $S_n=\sum_{k=1}^n [ g_{\theta_k}(X_k)- \pi_{\theta_k}(g_{\theta_k}(X_k))]$, where 
$\pi_{\theta_k}(g_{\theta_k}(X_k))=\int_{\BX} g_{\theta_k}(x)) \pi_{\theta_k}(x) dx$. 
Let $S_n'=\sum_{k=1}^n [u_{\theta_k}-P_{\theta_k} u_{\theta_k}]$, where 
 $u_{\theta_k}$ is the solution to the Poisson equation 
\[
 u_{\theta_k}-P_{\theta_k} u_{\theta_k} = g_{\theta_k}(X_k)-\pi_{\theta_k}(g_{\theta_k}(X_k)).
\]
Further, we decompose $S_n$ into three terms, $S_n=S_n^{(1)}+S_n^{(2)}+S_n^{(3)}$, where 
\begin{align*}
 S_n^{(1)}=&\sum_{k=1}^n\left[ u_{\theta_{k-1}}(X_k)-P_{\theta_{k-1}}u_{\theta_{k-1}}(X_{k-1}) \right ],\\
 S_n^{(2)}=&\sum_{k=1}^n\left[ u_{\theta_k}(X_k)-u_{\theta_{k-1}}(X_k) \right ], \\
 S_n^{(3)}=&P_{\theta_0}u_{\theta_0}(X_0)-P_{\theta_n}u_{\theta_n}(X_n). 
\end{align*}

By Lemma \ref{lem60}, for all $\theta$ and $X$, there exists a constant $c$ such that 
 \[
 |u_\theta(X)| \leq c, \quad \mbox{and} \quad 
|P_\theta u_\theta(X)| \leq c.
\]
Let $p>2$, and $(1+\tau')/2\leq \iota'<\iota$ (where $\tau'$ is defined in ($A_4$)). 
Thus, there exists a constant $c$ such that
\[
 E\left\{|u_{\theta_{k-1}}(X_k)+P_{\theta_{k-1}}u_{\theta_{k-1}}(X_{k-1})|^p  \right \}\leq c .
\]
 Since $E\left[ u_{\theta_{k-1}}(X_k)-P_{\theta_{k-1}}u_{\theta_{k-1}}(X_{k-1})|\mathcal F_{k-1}\right] 
 =P_{\theta_{k-1}}u_{\theta_{k-1}}(X_{k-1})-P_{\theta_{k-1}}u_{\theta_{k-1}}(X_{k-1})=0$, 
$\{S_n^{(1)} \}$ is a martingale with the increments upper bounded in $L^p$.
Hence, by Burkholder inequality (Hall and Heyde, 1980; Theorem 2.10) and Minkowski's inequality,
 there exists a constant $c$ and $c'$ such that 
\begin{align*}
 E\left\{ |S_n^{(1)}|^p\right\} & \leq c E\left\{\left(\sum_{k=1}^n|u_{\theta_{k-1}}(X_k)-P_{\theta_{k-1}}u_{\theta_{k-1}}(X_{k-1})|^2\right)^{p/2}
  \right\}\\
 &\leq c \left\{\sum_{k=1}^n \left( E\left[ |u_{\theta_{k-1}}(X_k)-P_{\theta_{k-1}}u_{\theta_{k-1}}(X_{k-1})|^p\right]
  \right)^{2/p} \right\}^{p/2} \\
 &\leq c' n^{p/2} .
\end{align*}

Now we consider $S_n^{(2)}$. 
By Lemma \ref{lem60}, the fact $\|\theta_k-\theta_{k-1}\|= \gamma_k \| H(\theta_{k-1},X_k) \|$ and $(A_3)$-(ii),   
\[
|S_n^{(2)}|=\left|\sum_{k=1}^n \{u_{\theta_k}(X_k)-u_{\theta_{k-1}}(X_k)\} \right| 
 \leq c \sum_{k=1}^n \|\theta_k-\theta_{k-1}\|^{\iota'} 
\leq c'\sum_{k=1}^n\gamma_k^{\iota'},
\] 
Hence, $E(|S_n^{(2)}|^p)\leq c'^{p} (\sum_{k=1}^n\gamma_k^{\iota'})^p$.
Also, the third term is also bounded by some constant $c$, $E(|S_n^{(3)}|^p)<c$.

Hence, by Minkowski's inequality, Markov's inequality,
we can conclude
\begin{equation}\label{equm2}
 P\{n^{-1}|S_n|\geq\delta\}\leq C\delta^{-p}\left\{n^{-p/2}+(n^{-1}\sum_{k=1}^n\gamma_k^{\iota'})^p+n^{-p}\right\},
\end{equation}
where $C$ denotes a constant.
By (\ref{equm2}) and the Borel-Cantelli lemma, we have 
\[
 P\{\sup_{n\geq 1} n^{-1}|S_n|\geq\delta  \}\leq
 C\delta^{-p}\sum_{n\geq 1}\left\{n^{-p/2}+n^{-p/2}(n^{-1/2}\sum_{k=1}^n\gamma_k^{\iota'})^p+n^{-p}\right\}.
\]
Then, the SLLN is concluded with Kronecker's lemma, condition (\ref{coneq1}) and the condition $p>2$.
\end{proof}

\paragraph{Proof of Lemma \ref{lem51}} 
\begin{itemize}
  \item[(i)] Define 
\begin{equation} \label{noisedecomeq}
\begin{split} 
    e_{k+1} &= u_{\theta_k}(X_{k+1})-P_{\theta_k} u_{\theta_k}(X_k) , \\
   \nu_{k+1}&=\big[ P_{\theta_{k+1}} u_{\theta_{k+1}}(X_{k+1})
    -P_{\theta_k} u_{\theta_{k}}(X_{k+1}) \big ]+ \frac{\gamma_{k+2}-\gamma_{k+1}}{\gamma_{k+1}}
     P_{\theta_{k+1}} u_{\theta_{k+1}}(X_{k+1}), \\
\tilde{\varsigma}_{k+1}&=\gamma_{k+1} P_{\theta_k} u_{\theta_k}(X_k), \\
\varsigma_{k+1}&=\frac{1}{\gamma_{k+1}} (\tilde{\varsigma}_{k+1}-\tilde{\varsigma}_{k+2}), \\
\end{split}
\end{equation}
 where $u_{\cdot}(\cdot)$ is the solution of the Poisson equation (Refer to Lemma \ref{lem1}).
 It is easy to verify that 
 $H(\theta_{k},X_{k+1})-h(\theta_k)=e_{k+1}+\nu_{k+1}+\varsigma_{k+1}$ holds.

\item[(ii)] By (\ref{noisedecomeq}), we have
\begin{equation} \label{Gammaeq1}
E(e_{k+1}|\mF_k)=E(u_{\theta_k}(X_{k+1})|\mF_k)-P_{\theta_k}u_{\theta_k}(X_k)=0,
\end{equation}
 Hence, $\{e_k\}$ forms a martingale difference sequence.
 Following from Lemma \ref{lem1}-($B_2$), we have 
 \begin{equation} \label{Gammaeq2}
  \sup_{k\geq 0} E(\|e_{k+1}\|^{\alpha} |\mF_k) 1_{\{ \|\theta_k-\theta_*\| \leq \rho \}}<\infty.
 \end{equation}
 This concludes part (ii).

\item[(iii)]
 By (\ref{noisedecomeq}), we have
\begin{equation} \label{lem3proofeq1}
 E(e_{k+1} e_{k+1}^T|\mF_{k})=E\left[ u_{\theta_k}(X_{k+1})u_{\theta_k}(X_{k+1})^T | \mF_{k} \right] -
  P_{\theta_k} u_{\theta_k}(X_k) P_{\theta_k} u_{\theta_k}(X_k)^T
 \stackrel{\triangle}{=} l(\theta_k,X_k).
\end{equation}
By ($B_2$) and ($A_3$)-(i), there exist constants $c_1$, $c_2$, $c_3$ and $M$ 
such that
\begin{align*}
 \|l(\theta_k,X_k)\|&\leq E(\|u_{\theta_k}(X_{k+1})u_{\theta_k}(X_{k+1})^T\||\mathcal F_k)
                 +\|P_{\theta_k}u_{\theta_k}(X_k)P_{\theta_k}u_{\theta_k}(X_k)^T\|<c_1
\end{align*}
For any $\theta_k$, $\theta_k'\in\Theta$,
\begin{align}
 \|l(\theta_k,X_k)-l(\theta'_k,X_k)\|&\leq E(\|u_{\theta_k}(X_{k+1})u_{\theta_k}(X_{k+1})^T-u_{\theta_k'}(X_{k+1})u_{\theta_k'}(X_{k+1})^T\||\mathcal F_k)\nonumber\\
&+\|P_{\theta_k}u_{\theta_k}(X_k)P_{\theta_k}u_{\theta_k}(X_k)^T
 -P_{\theta_k'}u_{\theta_k'}(X_k)P_{\theta_k'}u_{\theta_k'},X_k)^T\label{eqA}\|.
\end{align}
By lemma \ref{lem1} ($B_2$)-($ii$), we have, for any $\eta\in(0,1)$  
\begin{align*}
 &\|P_{\theta_k}u_{\theta_k}(X_k)P_{\theta_k}u_{\theta_k}(X_k)^T-P_{\theta_k'}u_{\theta_k'}(X_k)P_{\theta_k'}u_{\theta_k'}(X_k)^T\|\\
\leq & \|(P_{\theta_k}u_{\theta_k}(X_k)-P_{\theta_k'}u_{\theta_k'}(X_k))P_{\theta_k}u_{\theta_k}(X_k)^T\|
 +\|P_{\theta_k'}u_{\theta_k'}(X_k)(P_{\theta_k}u_{\theta_k}(X_k)^T-P_{\theta_k'}u_{\theta_k'}(X_k)^T)\|\\
\leq& c_2 \|\theta_k-\theta_k'\|^\eta,
\end{align*}
 and 
\begin{align*}
 &E(\|u_{\theta_k}(X_{k+1})u_{\theta_k}(X_{k+1})^T-u_{\theta_k'}(X_{k+1})u_{\theta_k'}(X_{k+1})^T\||\mathcal F_k)\\
\leq & E(\|(u_{\theta_k}(X_{k+1})-u_{\theta_k'}(X_{k+1}))u_{\theta_k}(X_{k+1})^T\||\mathcal F_k)
 +E(\|u_{\theta_k'}(X_{k+1})(u_{\theta_k}(X_{k+1})^T-u_{\theta_k'}(X_{k+1})^T)\||\mathcal F_k)\\
\leq & c_3 \|\theta_k-\theta_k'\|^\eta .
\end{align*}
Plug into equation (\ref{eqA}), we have $\|l(\theta_k,X_k)-l(\theta_k',X_k)\|\leq M \|\theta_k-\theta_k'\|^\eta$ for any 
 $\theta_k, \theta_k'\in \Theta$, where $M$ is a constant.  

Let $\iota=\eta\in((\tau'+1)/2,1)$, 
then the conditions of Lemma \ref{slln} hold and thus 
\begin{equation}
\frac{1}{n}\sum_{k=1}^n\left[ l(\theta_k,X_k)- \pi_{\theta_k}(l(\theta_k,X)) \right] \rightarrow 0, \quad a.s.
\end{equation}
 where $\pi_{\theta_k}(l(\theta_k,X))=\int_{\BX} l(\theta_k,x) \pi_{\theta_k}(x) dx$.

On the other hand, we have 
\begin{align*}
 & \left \|\pi_{\theta_k}(l(\theta_k,X))-\pi_{\theta_*}(l(\theta_*,X)) \right \|\\
\leq& \left \| \pi_{\theta_k}\left( l(\theta_k,X)- l(\theta_*,X) \right) \right \|+
 \left \| \pi_{\theta_k} (l(\theta_*,X)) - \pi_{\theta_*}( l(\theta_*,X))  \right \|\\
\leq& M \left \|\theta_k-\theta_* \right \|^\eta 
 +\left \| \pi_{\theta_k}(l(\theta_*,X))- \pi_{\theta_*}( l(\theta_*,X)) \right \|.
\end{align*}
Given $\theta_k\rightarrow\theta_*$ a.s.,
 the first term goes to 0 almost surely as $k \to \infty$. 
 By condition ($A_3$), which implies the conditions of proposition 1.3.6  of Atchad$\acute{\rm e}$ {\it et al.} (2011) holds, 
 therefore  $\pi_{\theta_k}(l(\theta_*,X))- \pi_{\theta_*}( l(\theta_*,X)) \to 0$ almost surely. 
Thus, $\|\int_{\BX} l(\theta_k,x)d\pi_{\theta_k}(x)-\int_{\BX} l(\theta_*,x)d\pi_{\theta_*}(x)\|\rightarrow0$ almost surely and 
\begin{equation}\label{Gamma}
 \frac{1}{n}\sum_{k=1}^nl(\theta_k,X_k)\rightarrow\int_{\BX} l(\theta_*,x)d\pi_{\theta_*}(x)=\Gamma, \quad a.s.
\end{equation}
for some positive definite matrix $\Gamma$.
This concludes  part (iii).

\item[(iv)] By condition $(A_4)$, we have
 \[
 \frac{\gamma_{k+2}-\gamma_{k+1}}{\gamma_{k+1}}= O(\gamma_{k+2}^{\tau}),
 \]
 for some value $\tau \in [1,2)$. 
 By (\ref{noisedecomeq}) and (\ref{boundeq1}), there exists a
 constant $c_1$ such that the following inequality holds,
 \[
 \|\nu_{k+1}\| \leq c_1 \|\theta_{k+1}-\theta_{k}\|+ O(\gamma_{k+2}^{\tau})
  =c_1 \|\gamma_{k+1} H(\theta_k, X_{k+1}) \| + O(\gamma_{k+2}^{\tau}),
 \]
 which implies, by (\ref{Drieq21}), that there exists a constant $c_2$ such that
 \begin{equation} \label{nunormeq}
 \|\nu_{k+1}\| \leq c_2 \gamma_{k+1}.
 \end{equation}
therefore,
 \[
 E(\|\nu_k\|^2/\gamma_k) 1_{\{ \|\theta_k-\theta_*\| \leq \rho \}} \to 0.
 \]
 This concludes part (iv).
 
 \item[(v)]
 A straightforward calculation shows that 
 \[
 \gamma_{k+1} \varsigma_{k+1}=\tilde{\varsigma}_{k+1}-\tilde{\varsigma}_{k+2}
  =\gamma_{k+1} P_{\theta_k} u_{\theta_k}(X_k)-\gamma_{k+2} P_{\theta_{k+1}} u_{\theta_{k+1}}(X_{k+1}),
\]
By $(B_2)$, $E \left[ \| P_{\theta_k} u_{\theta_k}(X_k)\| \right]$ 
is uniformly bounded with respect to $k$. Therefore, (v) holds.
\end{itemize}

\subsection{Proof of Theorem \ref{normalitytheorem}}

 To prove Theorem \ref{normalitytheorem}, we introduce Lemma \ref{lem52}, which a 
 combined restatement of Theorem D.6.4 (Meyn and Tweedie, 2009; p.563) and 
 Theorem 1 of Pelletier (1998).  

 \begin{lemma}\label{lem52} Consider a stochastic approximation algorithm of the form 
  \[
   Z_{k+1}=Z_k+\gamma_{k+1} h(Z_k)+ \gamma_{k+1}(\nu_{k+1}+e_{k+1}),
  \]
  where $\nu_{k+1}$ and $e_{k+1}$ are noise terms.  
  Assume that $\{\nu_k\}$ and $\{e_k\}$ satisfies (ii)-(iv) given in Lemma \ref{lem51},
  and the conditions $(A_2)$ and $(A_4)$ are satisfied.  
  On the set $\Lambda(z^*)=\{ Z_k \to z^*\}$, 
  \[
   \frac{Z_k-z^*}{\sqrt{\gamma_k}} \Longrightarrow \mN(0, \Sigma),
  \]
  with $\Longrightarrow$ denoting the weak convergence, $\mN$ the Gaussian distribution and
 \[
 \Sigma=\int_0^{\infty} e^{(F'+\zeta I)t} \Gamma e^{(F+\zeta I)t} dt,
 \]
 where $F$ is defined in $(A_2)$, $\zeta$ is defined in (\ref{Peleq2}), and
 $\Gamma$ is defined in Lemma \ref{lem51}.
 \end{lemma}

 \paragraph{Proof of Theorem \ref{normalitytheorem}}
  Rewrite the SAMCMC algorithm in the form 
 \begin{equation} \label{itertheta}
\theta_{k+1}-\theta_*=(\theta_{k}-\theta_*)+\gamma_{k+1} h(\theta_k)+ \gamma_{k+1}\xi_{k+1}.
\end{equation}
To facilitate the theoretical analysis for the random process $\{ \theta_k \}$,
 we define a reduced random process $\{ \tilde{\theta}_k \}_{k \geq 0}$:    
\begin{equation} \label{redproc}
 \tilde{\theta}_k =\theta_k+ \tvarsigma_{k+1},
\end{equation}
 where $\tvarsigma_{k+1}$ is as defined  in equation (\ref{noisedecomeq}) in the proof of Lemma \ref{lem51}.
 Then, for the SAMCMC algorithm, we have
\begin{align} 
\tilde{\theta}_{k+1}-\theta_*&=(\tilde{\theta}_{k}-\theta_*)+\gamma_{k+1} h({\theta}_{k})+\gamma_{k+1}\xi_{k+1}+\tvarsigma_{k+2}-\tvarsigma_{k+1} \nonumber
\\&=(\tilde{\theta}_{k}-\theta_*)+\gamma_{k+1} h(\tilde{\theta}_{k})+\gamma_{k+1}(h(\theta_k)-h(\tilde\theta_k)+\xi_{k+1}-\varsigma_{k+1}) \nonumber 
\\&= (\tilde{\theta}_{k}-\theta_*)+\gamma_{k+1} h(\tilde{\theta}_{k})+
  \gamma_{k+1}(h(\theta_k)-h(\tilde\theta_k)+v_{k+1}+e_{k+1}) \nonumber
\\&=(\tilde{\theta}_{k}-\theta_*)+\gamma_{k+1} h(\tilde{\theta}_{k})+ 
\gamma_{k+1}(\tilde{\nu}_{k+1}+ e_{k+1}), \label{itertheta2}
\end{align}
 where $\tilde{\nu}_{k+1}= \nu_{k+1}+h(\theta_k)-h(\tilde{\theta}_{k})$, and  $\varsigma_{k+1}$,
 $\nu_{k+1}$ and $e_{k+1}$ are defined in  equation (\ref{noisedecomeq}) in the proof of Lemma \ref{lem51} as well.
 Since $h(\cdot)$ is H$\ddot{\rm o}$lder continuous on $\Theta$ (by the result $B_3$ of Lemma \ref{lem1}) and $\Theta$ is compact,  
 there exists a constant $M$ such that $\|h(\theta_k)-h(\tilde{\theta}_{k})\|
 \leq M\|\tilde{\theta}_{k}-\theta_{k}\|^{\eta}=M\|\tvarsigma_{k+1}\|^{\eta}$ for any $\eta \in (0.5,1)$. 
Thus, by (\ref{noisedecomeq}), 
 there exists a constant $c$ such that 
 \[
 E \left[ \|h(\theta_k)-h(\tilde{\theta}_{k})\|^2/\gamma_k\right] \leq c \gamma_{k+1}^{2 \eta-1}\frac{\gamma_{k+1}}{\gamma_k}\rightarrow 0 ,
 \]
 since  $\gamma_{k+1}^{2 \eta-1} \to 0$ and $\gamma_{k+1}/\gamma_k \to 1$ as $k \to \infty$.

 Therefore, $\tilde{\nu}_{k+1}=\nu_{k+1}+h(\theta_k)-h(\tilde{\theta}_{k})$ also satisfies 
 the property (iv) of Lemma \ref{lem51}.

 By Lemma \ref{lem51} and Lemma \ref{lem52}, we have 
 \[
   \frac{\tilde{\theta}_k-\theta_*}{\sqrt{\gamma_k}} \Longrightarrow \mN(0, \Sigma).
 \]
 By Lemma \ref{lem1} , $E \| P_{\theta_k} u_{\theta_k}(X_k) \|$ 
 is uniformly bounded with respect to $k$. Hence, 
 \begin{equation} \label{sluitskyeq}
 \frac{\tvarsigma_{k+1}}{\sqrt{\gamma_k}} \to 0, \quad \quad \mbox{in probability}.
 \end{equation} 
 It follows from Slutsky's theorem (see, e.g., Casella and Berger, 2002), 
 \[
   \frac{\theta_k-\theta_*}{\sqrt{\gamma_k}} \Longrightarrow \mN(0, \Sigma),
 \]
 which concludes Theorem \ref{normalitytheorem}.

 \subsection{Proof of Theorem \ref{efftheorem}} 
 
 \paragraph{Proof of Theorem \ref{efftheorem}} 
 Let $\bx=(x^{(1)}, \ldots, x^{(\kappa)})$ denote the samples drawn at an iteration of population SAMCMC. 
 Let $\bP(\bx,\by)$ and $P(x,y)$ denote the Markovian transition kernels used in the population and 
 single-chain SAMCMC algorithms, respectively. Let $\bH(\theta,\bx)$ and $H(\theta,x)$ be the parameter updating 
 function associated with the population and single-chain SAMCMC algorithms, respectively.  
 Let $\bu=\sum_{n \geq 0} (\bP^n \bH-h)$ be a solution of Poisson equation $\bu-\bP \bu=\bH-h$,
 and let $u=\sum_{n \geq 0} (P^n H-h)$ be a solution of Poisson equation $u-P u=H-h$. 
 Since 
 \[
 \bH(\theta,\bx)= \frac{1}{\kappa} \sum_{i=1}^{\kappa} H(\theta,x^{(i)}),
 \]
 we have $\bu_{\theta}(\bx)=\frac{1}{\kappa} \sum_{i=1}^{\kappa} u_{\theta}(x^{(i)})$. 
 By (\ref{noisedecomeq}), we further have 
 \[
  \be_{t+1}=\frac{1}{\kappa} \sum_{i=1}^{\kappa} e_{t+1}^{(i)}.
 \]
 Since $x_{t+1}^{(1)}, \ldots, x_{t+1}^{(\kappa)}$ are mutually independent conditional on 
 $\mF_t$, $e_{t+1}^{(1)}, \ldots, e_{t+1}^{(\kappa)}$ are also independent conditional on $\mF_t$ and thus 
 \[
 \bGamma=\Gamma/\kappa,
 \]
 which, by Theorem \ref{normalitytheorem}, further implies
 \[
 \Sigma_p=\Sigma_s/\kappa,
 \]
 where $\Sigma_p$ and $\Sigma_s$ denote the limiting covariance matrices
 of population SAMCMC and single-chain SAMCMC algorithms, respectively. 
 Therefore,  
 $(\theta_t^p-\theta_*)/\sqrt{\gamma_t}$ and $(\theta_{\kappa t}^s -\theta_*)/\sqrt{\kappa\gamma_{\kappa t}}$
 both converge in distribution to $N(0,\Sigma_p)$.
 By condition $(A_4)$, $\gamma_t/(\kappa \gamma_{\kappa t})= \kappa^{\beta-1}$, which concludes the proof.

{\centering \section*{\sc REFERENCES} }

\begin{description}

\item[] Aldous, D., Lov$\acute{\rm a}$sz, L., and Winkler, P. (1997). Mixing times for uniformly ergodic 
        Markov chains. {\it Stochastic Processes and Their Applications}, {\bf 71}, 165-185.

\item[] Andrieu, C. and Moulines, $\acute{\rm E}$ (2006). On the ergodicity properties of some adaptive MCMC algorithms. 
  {\it Annals of Applied Probability}, {\bf 16}, 1462-1505.

\item[] Andrieu, C., Moulines, $\acute{\rm E}$, and Priouret, P. (2005). Stability of Stochastic Approximation 
    Under Verifiable Conditions. {\it SIAM Journal of Control and Optimization}, {\bf 44}, 283-312.

\item[] Atchad$\acute{\rm e}$, Y. and Fort, G. (2009). Limit theorems for some adaptive MCMC algorithms 
        with subgeometric kernels. {\it Bernoulli}, {\bf 16}, 116-154.

\item[] Atchad$\acute{\rm e}$, Y.,  Fort, G. Moulines, E. and Priouret, P. (2011) Adaptive Markov chain 
   Monte Carlo: Theory and methods. In {\it Bayesian Time Series Models}. Cambridge University Press, Oxford, UK.

\item[] Benveniste, A., M$\acute{\rm e}$tivier, M., and Priouret, P. (1990). {\it Adaptive Algorithms and Stochastic
        Approximations}. New York: Springer-Verlag.

\item[] Billingsley, P. (1986). {\it Probability and Measure} (2nd edition). New York: John Wiley \& Sons. 

\item[] Blum, J.R. (1954). Approximation Methods which Converge with Probability one. 
        {\it Ann. Math. Statist.} {\bf 25}, 382-386. 

\item[] Casella, G. and Berger, R.L. (2002). {\it Statistical Inference} (second edition). Duxbury Thomson Learning. 


\item[] Chauveau, D. and Diebolt, J. (2000). Stability properties for a product Markov chain. 
        Preprint No 06/2000, Universit$\acute{\rm e}$ Marne-la-Vall$\acute{\rm e}$e. 

\item[] Chen, H.F. (2002). {\it Stochastic Approximation and Its Applications}. Dordrecht: Kluwer Academic Publishers.

\item[] Cheon, S. and Liang, F. (2009). Bayesian phylogeny analysis via stochastic approximation Monte Carlo.
 {\it Molecular Phylogenetic \& Evolution}, {\bf 53}, 394-403.

\item[] Duan, G.-R. and Patton, R.J. (1998). A Note on Hurwitz Stability of Matrices. {\it Automatica}, {\bf 34}, 509-511.

\item[] Geman, S., and Geman, D. (1984). Stochastic relaxation, Gibbs distributions and the Bayesian restoration of images.
  {\it IEEE Transactions on Pattern Analysis and Machine Intelligence}, {\bf 6}, 721-741.

\item[] Geyer, C.J. (1991). Markov chain Monte Carlo maximum likelihood.
In {\it Computing Science and Statistics: Proceedings of the 23rd Symposium on
 the Interface } (ed. E.M. Keramigas), pp.153-163.

\item[] Gilks, W.R., Roberts, G.O., and George, E.I. (1994). Adaptive Direction Sampling, 
       {\it The Statistician}, {\bf 43}, 179-189.

 \item[] Gu, M.G. and Kong, F.H. (1998). A stochastic approximation algorithm with Markov
 chain Monte Carlo method for incomplete data estimation problems.
{\it Proc. Natl. Acad. Sci. USA}, {\bf 95} 7270-7274.

\item[] Haario, H., Saksman, E., and Tamminen, J. (2001). An adaptive Metropolis algorithm. 
       {\it Bernoulli}, {\bf 7}, 223-242.

\item[] Hall, P. and Heyde, C. C. (1980). {\it Martingale limit theory and its applications}, 
       Academic Press, New York, London.

\item[] Hastings, W.K. (1970). Monte Carlo sampling methods using Markov chain and their applications.
        {\it Biometrika}, {\bf 57}, 97-109.
 
 \item[] Liang, F. (2007). Continuous contour Monte Carlo for marginal density estimation with an application
        to a spatial statistical model. \JCGS, {\bf 16}, 608-632

\item[] Liang, F. (2009). Improving SAMC Using Smoothing Methods: Theory and Applications to Bayesian Model Selection Problems.
      {\it The Annals of Statistics}, {\bf 37}, 2626-2654.

 \item[] Liang, F. (2010). Trajectory averaging for stochastic approximation MCMC algorithms. 
         {\it Annals of Statistics}, {\bf 38}, 2823-2856.

\item[]  Liang, F., Liu, C. and Carroll, R. J. (2007) Stochastic approximation in Monte Carlo computation.
  \JASA, {\bf 102}, 305-320.

\item[] Liang, F., and Wong, W.H. (2000). Evolutionary Monte Carlo: Application to $C_p$ model sampling and change point problem.
        {\it Statistica Sinica}, {\bf 10}, 317-342.

\item[] Liang, F., and Wong, W.H. (2001). Real parameter evolutionary Monte Carlo with applications in Bayesian 
        mixture models. \JASA, {\bf 96}, 653-666.

\item[] Liang, F. and Zhang, J. (2009). Learning Bayesian Networks for Discrete Data. 
    {\it Computational Statistics $\&$Data Analysis}, {\bf 53}, 865-876.

\item[] Liu, J.S., Liang, F., and Wong, W.H. (2000). The use of multiple-try method and local optimization in 
   Metropolis sampling. {\it Journal of American Statistical Association}, {\bf 94}, 121-134.

\item[] Marinari, E., and Parisi, G. (1992). Simulated Tempering: A New Monte Carlo Scheme.
       {\it Europhys. Lett.}, {\bf 19}, 451-458.

\item[] Metropolis N., Rosenbluth A.W., Rosenbluth M.N., Teller A.H., and Teller E. (1953).
       Equation of state calculations by fast computing machines.
       {\it Journal of Chemical Physics}, {\bf 21}, 1087-1091.

\item[] Meyn, S. and Tweedie, R.L. (2009). {\it Markov Chains and Stochastic Stability} (second edition). 
        Cambridge University Press. 

\item[] Nummelin, E. (1984), {\it General Irreducible Markov Chains and Nonnegative Operators.} Cambridge: Cambridge University Press.

\item[] Pelletier, M. (1998). Weak convergence rates for stochastic approximation with application to multiple
        targets and simulated annealing. {\it Annals of Applied Probability}, {\bf 8}, 10-44.

 \item[] Robbins, H. and Monro, S. (1951). A Stochastic approximation method.
  {\it Ann. Math. Statist.}, {\bf 22} 400-407.

\item[] Roberts, G.O. and Rosenthal, J.S. (2007). Coupling and ergodicity of adaptive Markov chain
  Monte Carlo algorithms. {\it J. Appl. Prob.}, {\bf 44}, 458-475.

\item[] Roberts, G.O., and Rosenthal, J.S. (2009). Examples of adaptive MCMC. 
  \JCGS, {\bf 18}, 349-367.

\item[] Roberts, G.O., and Tweedie, R.L. (1996). Geometric Convergence and Central Limit Theorems
 for Multidimensional Hastings and metropolis Algorithms. {\it Biometrika}, {\bf 83}, 95-110.

\item[] Tadi$\acute{\mbox{c}}$, V. (1997). On the convergence of stochastic iterative algorithms and
 their applications to machine learning. A short version of this paper was published in
  {\it Proc. 36th Conf. on Decision $\&$ Control} 2281-2286. San Diego, USA.

 \item[] Younes, L. (1989). Parametric inference for imperfectly observed Gibbsian fields.
 {\it Probab. Theory Relat. Field}, {\bf 82} 625-645.

\item[] Younes, L. (1999). On the convergence of Markovian stochastic algorithms with rapidly 
        decreasing ergodicity rates. {\it Stochastics and Stochastics Reports}, {\bf 65}, 177-228.

\item[] Wang, F. and Landau, D.P. (2001). Efficient, multiple-range random walk algorithm to calculate the density of states.
       {\it Physical Review Letters}, {\bf 86}, 2050-2053.

\item[] Wong, W.H. and Liang, F. (1997). Dynamic weighting in Monte Carlo and optimization.
  \PNAS, {\bf 94}, 14220-14224.


\item[] Ziedan, I.E. (1972). Explicit solution of the Lyapunov-matrix equation. {\it IEEE Transactions on Automatic
        Control}, {\bf 17}, 379-381.
\end{description}

\end{document}


\thispagestyle{empty}

\title{Supplementary Material for ``Weak Convergence Rates of Population versus Single-Chain Stochastic Approximation MCMC Algorithms''}

\author{Qifan Song,  Mingqi Wu, Faming Liang
 \thanks{correspondence author: Faming Liang.
 Faming Liang is professor, 
 Department of Statistics, Texas A$\&$M University, College Station, TX 77843,
 Email: fliang@stat.tamu.edu;
 Qifan Song is a graduate student,
 Department of Statistics, Texas A$\&$M University, College Station, TX 77843.
 Mingqi Wu is statistical consultant, Shell Global Solutions (US) Inc., 
Shell Technology Center Houston,
3333 Highway 6 South, Houston, TX 77082-3101.
 }}

\maketitle

\section{Proof of Theorem 2.1} 
 
\begin{proof} Let $M=\sup_{\theta \in \Theta} \max\{ \|h(\theta)\|,|v(\theta)|\}$ and
 $\mV_{\varepsilon}=\{\theta: d(\theta, \mathcal L) \leq \varepsilon \}$.
 Applying Taylor's expansion formula (Folland, 1990), we have
\[
v(\theta_{t+1})=v(\theta_t)+\gamma_{n+1} v_h(\theta_{t+1})+R_{t+1}, \quad t\geq 0,
\]
which implies that
\[
\sum_{i=0}^t \gamma_{i+1} v_h(\theta_i)=v(\theta_{t+1})-v(\theta_0)-\sum_{i=0}^t R_{i+1}
 \geq -2 M-\sum_{i=0}^{t} R_{i+1}.
\]
Since $\sum_{i=0}^{t} R_{i+1}$ converges (owing to Lemma A.2),
 $\sum_{i=0}^t \gamma_{i+1} v_h(\theta_i)$ also converges. Furthermore,
\[
v(\theta_t)=v(\theta_0)+\sum_{i=0}^{t-1} \gamma_{i+1} v_h(\theta_i)+\sum_{i=0}^{t-1} R_{i+1}, \quad t \geq 0,
\]
$\{ v(\theta_t)\}_{t \geq 0}$ also converges. On the other hand, conditions $(A_1)$ and $(A_2)$ imply
$\varliminf_{t \rightarrow \infty} d(\theta_t,\mL)=0$. Otherwise, there exists $\varepsilon>0$ and $n_0$
 such that $d(\theta_t,\mL) \geq \varepsilon$, $t \geq n_0$; as $\sum_{t=1}^{\infty} \gamma_{t}=\infty$
 and $p=\sup\{ v_h(\theta): \theta \in \mV_{\varepsilon}^c \}<0$, it is obtained that
$\sum_{t=n_0}^{\infty} \gamma_{t+1}v_h(\theta_t) \leq p \sum_{t=1}^{\infty} \gamma_{t+1}=-\infty$.

 Suppose that $\varlimsup_{t \rightarrow \infty} d(\theta_t,\mL)>0$. Then, there exists $\varepsilon>0$
 such that $\varlimsup_{t \rightarrow \infty} d(\theta_t,\mL) \geq 2 \varepsilon$.
 Let $t_0=\inf\{t\geq 0: d(\theta_t,\mL) \geq 2 \varepsilon\}$, while
 $t_k'=\inf\{ t  \geq t_k: d(\theta_t,\mL) \leq \varepsilon\}$ and
 $t_{k+1}=\inf\{t \geq t_k': d(\theta_t,\mL) \geq 2 \varepsilon\}$, $k \geq 0$. Obviously,
 $t_k < t_{k'} < t_{k+1}$, $k \geq 0$, and
 \[
 d(\theta_{t_k},\mL)\geq 2 \varepsilon, \ d(\theta_{t_k'},\mL)\leq \varepsilon, \ \mbox{and} \
 d(\theta_t,\mL) \geq \varepsilon, \ t_k \leq t < t_k', \ k\geq 0.
\]
Let $q=\sup\{v_h(\theta): \theta \in \mV_{\varepsilon}^c\}$. Then
\[
q \sum_{k=0}^{\infty} \sum_{i=t_k}^{t_k'-1} \gamma_{i+1} \geq \sum_{k=0}^{\infty} \sum_{i=t_k}^{t_k'-1}
 \gamma_{i+1} v_h(\theta_i) \geq \sum_{t=0}^{\infty} \gamma_{t+1} v_h(\theta_t) > -\infty.
\]
Therefore, $\sum_{k=0}^{\infty} \sum_{i=t_k}^{t_k'-1} \gamma_{i+1}<\infty$, and consequently,
$\lim_{k \rightarrow \infty} \sum_{i=t_k}^{t_k'-1} \gamma_{i+1}=0$. Since
$\sum_{t=1}^{\infty} \gamma_t \xi_t$ converges (owing to Lemma A.2), we have
\[
\varepsilon \leq \|\theta_{t_k'}-\theta_{t_k} \| \leq M \sum_{i=t_k}^{t_k'-1} \gamma_{i+1}
 + \left \| \sum_{i=t_k}^{t_k'-1} \gamma_{i+1} \xi_{i+1} \right \| \longrightarrow 0, 
\]
 as $k \rightarrow \infty$.  This contradicts with our assumption $\varepsilon>0$. Hence,
  $\varlimsup_{t \rightarrow \infty} d(\theta_t,\mL)>0$ does not hold. Therefore,
 $\lim_{t \rightarrow \infty} d(\theta_t,\mL)=0$ almost surely.
\end{proof}

\section{Proofs of Theorems for Pop-SAMC}  

In order to study the convergence of the Pop-SAMC algorithm, we introduce an equivalent variation of the Pop-SAMC algorithm.
Without loss of generality, we assume that $E_1$, $\dots$ ,$E_{m_0}$ are nonempty subregions, 
 and $E_{m_0+1}, \ldots, E_{m}$ are all empty.

\begin{itemize}
\item[1.] (Population sampling) 
  The sampling step is the same as described in Section 3.2 of the main text.

\item[2'.] (Weight updating) Set 
\begin{equation} \label{eq63}
  \theta_{t+1}=\theta_t+\gamma_{t+1} \tilde\bH(\theta_t, \bx_{t+1}), 
\end{equation}
where $\tilde\bH(\theta_t,\bx_{t+1})=\sum_{i=1}^{\kappa} \tilde H(\theta_t, x_{t+1}^{(i)})/\kappa$, and 
 $\tilde H(\theta_t,x_{t+1}^{(i)})=\bz_{t+1}-\bpi-(I(x_{t+1}^{(i)}\in E_{m_0})-\pi_{m_0})\utwi{1}$. 
where $\bz_{t+1}$ and $\bpi$ are as specified in the SAMC algorithm, and $\utwi{1}$ denotes a vector of 1s.
\end{itemize}

This variational Pop-SAMC algorithm adds a constant vector 
 $-\gamma_{t+1} \sum_{i=1}^\kappa(I(x^{(i)}_{t+1}\in E_{m_0})-\pi_{m_0})\utwi{1}/\kappa$ to the
 estimate of $\btheta$ of the original algorithm and thus keeps $\theta^{(m_0)}$ unchanged, say $\theta^{(m_0)}_t\equiv 0$.
 Hence, below we
only need to prove that Theorem 3.1 and Theorem 3.2 are true
for this variational Pop-SAMC algorithm.

\subsection{Proof of Theorem 3.1}

 Since $E_{m_0+1}$, $\dots$, $E_m$ are empty, $\theta_{m_0+1}$, $\dots$, $\theta_m$ are auxiliary
 variable, which do not affect the updating of $(\theta_i)_{i=1}^{m_0-1}$ and sampling step at all.
 Therefore, we can view the algorithm as of it is only update $(\theta_1,\dots,\theta_{m_0-1})^T$ with function
 $(\tilde \bH^{(1)},\dots,\tilde \bH^{(m_0-1)})^T$. 
 Once we prove that for $i=1,\dots,m_0-1$, 
 \[\theta_t^{(i)}\rightarrow \log\left(\int_{E_i} \psi(x) dx\right)-\log(\pi_i+\nu)-
\log\left(\int_{E_{m_0}} \psi(x) dx\right)+\log(\pi_{m_0}+\nu),\]
 almost surely, then it  is trivial to see that $\theta_t^{(i)}\rightarrow-\infty$ for $i>m_0$.
 (Because $\sum_{j=1}^{t}I(x_j^{(k)}\in E_{m_0})/t\rightarrow \pi_{m_0}+\nu$, and $\theta_t^{(i)}=-t\pi_i-\sum_{k=1}^\kappa
 \sum_{j=1}^{t}I(x_j^{(k)}\in E_{m_0})/\kappa+t\pi_{m_0}$
 for any $i>m_0$.)

 To prove the convergence of $\theta_t^{(i)}$ for $i< m_0$,
 it follows from Theorem 2.1 that we only need to verify that 
 Pop-SAMC satisfies the conditions $(A_1)$, $(A_3)$ and $(A_4)$.
 This is done as follows.
 
\begin{itemize}
\item[$(A_1)$] This condition can be verified as in Liang {\it et al.} (2007). Since a part of the proof
 will be used in proving Theorem 3.2, we re-produce the proof below.
Since the invariant distribution of the kernel $\bP_{\theta_t}(\bx,\by)$ is $f_{\theta_t}(\bx)$,
 for any fixed value of $\theta$, we have 
\begin{align} 
E(\tilde \bH^{(i)}(\theta,\bx))
 &= \frac{ \int_{E_i} \psi(x) d x/e^{\theta_{i}}-\int_{E_{m_0}} \psi(x) d x/e^{\theta_{m_0}}}{ \sum_{k=1}^m
 [\int_{E_k} \psi(x) d x/e^{\theta_{k}}]}-\pi_i+\pi_{m_0} \nonumber\\
 &=\frac{S_i-S_{m_0}}{S}-\pi_i+\pi_{m_0}, \quad i=1, \ldots, m_0-1,\label{app22}
\end{align}
where $\tilde\bH^{(i)}(\theta, \bx)$ denotes the $i$th component of $\tilde\bH(\theta,\bx)$,
 $S_i= \int_{E_i} \psi(x) d x/e^{\theta_{i}}$ and $S=\sum_{k=1}^{m_0} S_k$.
Thus,
\[
h(\theta)=\int_{\mX} H(\theta,\bx) f(d \bx)=
 \left(\frac{S_1}{S}-\pi_1, \ldots, \frac{S_{m_0-1}}{S}-\pi_{m_0-1} \right)^T-\frac{S_{m_0}}{S}+\pi_{m_0}.
\]

 It follows from (\ref{app22}) that $h(\theta)$ is a continuous function of $\theta$.
Let 
\[
v(\theta)=\frac{1}{2} \sum_{k=1}^{m_0} \left(\frac{S_k}{S}- \pi_k \right)^2,
 \]
 which, as shown below, has continuous partial derivatives of the first order.

Solving the system of equations formed by (\ref{app22}), we have
\[
\mL=\left\{(\theta_1,\ldots,\theta_{m_0-1}): \theta_{i}=\mbox{Const}+\log \Big(\int_{E_i} \psi(x) d x\Big)-\log(\pi_i+\nu), i=1, \ldots, m;
 \theta \in \Theta \right \},
\]
where $\mbox{Const}=\log(\pi_{m_0}+\nu)-\log\int_{E_{m_0}}\psi(x) dx$.
It is obvious that $\mL$ is nonempty and $v(\theta)=0$ for every $\theta \in \mL$.

To verify the conditions related to $\nabla v(\theta)$, we have the following calculations:
\begin{equation} \label{der0}
\begin{split}
  \frac{\partial S}{\partial \theta_i}= \frac{\partial S_i}{\partial \theta_i}  = -S_i,
 &  \quad \quad  \frac{\partial S_i}{\partial \theta_j} =  \frac{\partial S_j}{\partial \theta_i}= 0, \\
 \frac{\partial \big( \frac{S_i}{S} \big)}{\partial \theta_i} =-\frac{S_i}{S}(1-\frac{S_i}{S}),
 & \quad \quad \frac{\partial \big( \frac{S_i}{S} \big)}{\partial \theta_j} =
  \frac{\partial \big( \frac{S_j}{S} \big)}{\partial \theta_i}=\frac{S_i S_j}{S^2},  \\
\end{split}
\end{equation}
for $i, j=1, \ldots,m_0-1$ and $i \ne j $.
\begin{equation} \label{der1}
\begin{split}
\frac{\partial v(\theta)}{\partial \theta_i} & = \frac{1}{2} \sum_{k=1}^{m_0}
 \frac{\partial (\frac{S_k}{S}- \pi_k)^2}{\partial \theta_i} \\
 & = \sum_{j=1}^{m_0} (\frac{S_j}{S}-\pi_j) \frac{S_i S_j}{S^2}- (\frac{S_i}{S}-\pi_i) \frac{S_i}{S} \\
 &= \mu_{\eta^*} \frac{S_i}{S}-(\frac{S_i}{S}-\pi_i) \frac{S_i}{S}, \\
\end{split}
\end{equation}
for $i=1, \ldots, m_0-1$, where $\mu_{\eta^*}= \sum_{j=1}^{m_0} (\frac{S_j}{S}-\pi_j) \frac{S_j}{S}$.
 Thus, 
\begin{equation} \label{negder2}
\begin{split}
v_h(\theta)&=\langle \nabla v(\theta), h(\theta) \rangle \\
  & = \mu_{\eta^*} \sum_{i=1}^{m_0-1} (\frac{S_i}{S}-\pi_i) \frac{S_i}{S} -
 \sum_{i=1}^{m_0-1} (\frac{S_i}{S}-\pi_i)^2 \frac{S_i}{S} 
 -\sum_{i=1}^{m_0-1}\left(\mu_{\eta^*} \frac{S_i}{S}-(\frac{S_i}{S}-\pi_i) \frac{S_i}{S}\right)\left(\frac{S_{m_0}}{S}-\pi_{m_0}\right) \\
 &= -\big\{ \sum_{i=1}^{m_0} (\frac{S_i}{S}-\pi_i)^2 \frac{S_i}{S}- \mu_{\eta^*}^2 \} \\
 &= -\sigma_{\eta^*}^2 \leq 0, \\
\end{split}
\end{equation}
where $\sigma_{\eta^*}^2$ denotes the variance of the discrete distribution defined in the following table,

\vspace{0.5mm}

\begin{center}
\begin{tabular}{c|c|c|c} \hline
State $(\eta^*)$ & $\frac{S_1}{S}-\pi_1$ & $\cdots$ & $\frac{S_{m_0}}{S}-\pi_m$ \\ \hline
Prob.  & $\frac{S_1}{S}$        & $\cdots$ & $\frac{S_{m_0}}{S}$  \\ \hline
\end{tabular}
\end{center}

\vspace{0.5mm}

If $\theta \in \mL$, $v_h(\theta)=0$. Otherwise, $v_h(\theta)<0$ and for any compact set 
 $\mK \subset \mL^c$, $\sup_{\theta \in \mL} v_h(\theta)<0$.

 \item[$(A_3)$]
 Let $\bx_{t+1}=(x_{t+1}^{(1)},\ldots, x_{t+1}^{(\kappa)})$, 
  which is a sample produced by $\kappa$ independent Markov chains
 on the product space $\BX=\mX \times \cdots \times \mX$ with the transition kernel
\[
\bP_{\theta_t}(\bx,\by)=P_{\theta_t}(x^{(1)},y^{(1)})P_{\theta_t}(x^{(2)},y^{(2)})\cdots
 P_{\theta_t}(x^{(\kappa)},y^{(\kappa)}),
\]
where $P_{\theta_t}(x,y)$ denotes a one-step MH kernel at a given value of $\theta_t$.
 Under the assumptions that both $\Theta$ and $\mX$ 
 are compact and the proposal distribution is local positive, it has been shown in Liang {\it et al.} (2007)
 that $P_{\theta}(x,y)$ satisfies the drift condition $(A_3)$. In what follows, we will show that 
 $\bP_{\theta}(\bx,\by)$ also satisfies $(A_3)$, given that $P_{\theta}(x,y)$ satisfies $(A_3)$.
 
 To simplify notations, in what follows we will drop the subscript $t$, denoting $x_t$ by $x$, 
 $\bx_t$ by $\bx$, and
  $\theta_t=(\theta_{t1},\ldots,\theta_{tm})$ by $\theta=(\theta_1,\ldots,\theta_m)$.
 Roberts and Tweedie (1996) (Theorem 2) show that if the target distribution is bounded away from 0 and $\infty$ 
 on every compact set of $\mX$, then 
 the MH chain with a proposal distribution satisfying the local positive condition
 is irreducible and aperiodic, and every nonempty compact set is small.
 It follows from this result that $P_{\theta}(x,y)$
 is irreducible and aperiodic, and thus $\bP_{\theta}(x,y)$ is
 also irreducible and aperiodic.

If $\mX$ is compact, and furthermore $f(x)$ is bounded away from 0 and $\infty$,
 by equation (20) of the main text, $f_{\theta}(x)$ is uniformly bounded away from 0 and $\infty$ since $\Theta$ is compact.
 By Roberts and Tweedie's arguments, these imply that
 $\mX$ is a small set and the minorization condition uniformly holds on $\mX$ for all kernel
 $P_{\theta}(x,y)$, $\theta\in\Theta$; i.e., there exist
 a constant $\delta$ and a probability measure $\nu'(\cdot)$ such that
\[
P_{\theta}(x,A) \geq \delta' \nu'(A), \quad \forall x \in \mX, \ \forall A \in \mB_{\mX}.
\]
Therefore, 
\[
\bP_{\theta}(\bx,\bA) \geq \delta \nu(\bA), \quad \forall \bx \in \BX,
 \ \forall \bA \in \mB_{\BX},
\]
where $\bA=A_1 \times A_2 \times \ldots \times A_{\kappa}$,  $\delta=(\delta')^{\kappa}$, and 
 $\nu(\bA)=\nu'(A_1) \times \nu'(A_2) \times \ldots \times \nu'(A_{\kappa})$.
Hence, ($A_3$-i) is satisfied.

For Pop-SAMC, we have $\bH(\theta,\bx)=\sum_{i=1}^{\kappa} H(\theta,x^{(i)})/\kappa$. 
Since each component of $\bH(\theta,\bx)$ takes a value between 0 and 1,  there exists
 a constant $c_1=\sqrt{m}$ such that for any $\theta \in \Theta$ and all $\bx \in \BX$,
\begin{equation} \label{proofeq1}
\|\bH(\theta,\bx)\| \leq c_1.
\end{equation}
Also, $\bH(\theta,\bx)$ does not depend on $\theta$ for a given sample $\bx$.
Hence, $\bH(\theta,\bx)-\bH(\theta',\bx)=0$ for all $(\theta,\theta')\in \Theta \times \Theta$, and
the following condition holds,
\begin{equation} \label{proofeq2}
 \|\bH(\theta,\bx)-\bH(\theta',\bx)\| \leq c_1 \|\theta-\theta'\|,
\end{equation}
for all $(\theta,\theta')\in \Theta \times \Theta$.
Equations (\ref{proofeq1}) and (\ref{proofeq2}) imply that ($A_3$-ii) is satisfied.

 In Liang {\it et al.} (2007), it has been shown for the single-chain MH kernel that 
 there exists a constant $c_2$ such that 
\begin{equation}  \label{proofeq4}
 |P_{\theta}(x,A)-P_{\theta'}(x,A)| \leq  c_2 \|\theta-\theta'\|,
\end{equation}
 for any measurable set $A \subset \mX$. Therefore, 
 there exists a constant $c_3$ such that
\[
\begin{split}
& \left|\bP_{\theta}(\bx,\bA)-\bP_{\theta'}(\bx,\bA)\right|= \big| \int_{A_1} \cdots \int_{A_{\kappa}}
\big[ P_{\theta}(x^{(1)},y^{(1)}) P_{\theta}(x^{(2)},y^{(2)}) \cdots P_{\theta}(x^{(\kappa)},y^{(\kappa)})-\\
 &  P_{\theta'}(x^{(1)},y^{(1)}) P_{\theta'}(x^{(2)},y^{(2)}) \cdots P_{\theta'}(x^{(\kappa)},y^{(\kappa)})
 \big] dy^{(1)} \cdots d y^{(\kappa)} \big| \\
 &\leq  \int_{A_1} \int_{\mX} \cdots \int_{\mX} \left|
  P_{\theta}(x^{(1)},y^{(1)})-P_{\theta'}(x^{(1)},y^{(1)})\right| P_{\theta}(x^{(2)},y^{(2)})
  \cdots P_{\theta}(x^{(\kappa)},y^{(\kappa)}) d y^{(1)} \cdots d y^{(\kappa)}  \\
 &  + \int_{\mX} \int_{A_2} \int_{\mX} \cdots \int_{\mX} P_{\theta'}(x^{(1)},y^{(1)})
 \left| P_{\theta}(x^{(2)},y^{(2)})- P_{\theta'}(x^{(2)},y^{(2)}) \right|
     P_{\theta}(x^{(3)},y^{(3)}) \cdots P_{\theta}(x^{(\kappa)},y^{(\kappa)}) \\
 & \ \ dy^{(1)}\cdots dy^{(\kappa)}
  +\cdots \\
 & + \int_{\mX} \cdots \int_{\mX} \int_{A_{\kappa}}
  P_{\theta'}(x^{(1)},y^{(1)})\cdots  P_{\theta'}(x^{(\kappa-1)},y^{(\kappa-1)}) \left|
   P_{\theta}(x^{(\kappa)},y^{(\kappa)})-
    P_{\theta'}(x^{(\kappa)},y^{(\kappa)})\right| dy^{(1)} \cdots d y^{(\kappa)}  \\
& \leq c_3 \|\theta-\theta'\|, \\
\end{split}
\]
which implies  $A_3$-(iii) is satisfied.

 \item[$(A_4)$]  $\{\gamma_t\}$  automatically satisfies the condition. 
\end{itemize}

\subsection{Proof of Theorem 3.2.}
 
 Following from Theorem 2.2 and Theorem 3.1, this 
 theorem can be proved by verifying that SAMC and Pop-SAMC satisfy $(A_2)$.
 To verify $(A_2)$, we first show that $h(\theta)$ has bounded first and second derivatives.
 Continuing the calculation
 in (\ref{der0}), we have
\begin{equation} \label{partialder}
\frac{\partial^2 ( \frac{S_i}{S} )}{\partial (\theta^{(i)})^2 } =\frac{S_i}{S} (1-\frac{S_i}{S})
 (1-\frac{2 S_i}{S}), \quad
\frac{\partial^2 ( \frac{S_i}{S} )}{\partial \theta^{(j)} \partial \theta^{(i)}}=
 -\frac{S_iS_j}{S^2} (1-\frac{2S_i}{S}),
\end{equation}
 where $S$ and $S_i$ are as defined in (\ref{app22}).
 This implies that the first and second derivatives of $h(\theta)$ are uniformly bounded by noting
 the inequality $0< \frac{S_i}{S}<1$. Hence, $h(\theta)$ is differentiable and its 
 derivative is Lipschitz continuous.

 Let $F=\partial h(\theta)/\partial \theta$. From (\ref{der0}) and (\ref{partialder}), we have
 $F=(\utwi{1}\utwi{1}^T+I)F_0$, where
 \[
 F_0=\begin{pmatrix}
 -\frac{S_1}{S}(1-\frac{S_1}{S}) & \frac{S_1 S_2}{S^2} & \cdots & \frac{S_1 S_{m_0-1}}{S^2} \\
 \frac{S_2S_1}{S^2} & -\frac{S_2}{S}(1-\frac{S_2}{S}) & \cdots & \frac{S_2S_{m_0-1}}{S^2} \\
    \vdots            &   \ddots                        & \vdots &  \vdots  \\
  \frac{S_{m_0-1}S_1}{S^2}&   \cdots   & \cdots & -\frac{S_{m_0-1}}{S}(1-\frac{S_{m_0-1}}{S}) \\
\end{pmatrix}.
\]
Thus, for any nonzero vector $\bz=(z_1,\ldots,z_{m_0-1})^T$, 
 \begin{equation} \label{SAMCCLTeqFull}
 \bz^T F_0 \bz =-\Big[ \sum_{i=1}^{m_0} z_i^2 \frac{S_i}{S} -\left(\sum_{i=1}^{m_0} z_i \frac{S_i}{S}
 \right)^2 \Big] =- \var(Z) < 0,
 \end{equation}
 where $z_{m_0}=0$, and $\var(Z)$ denotes the variance of the discrete
 distribution defined by the following table (note that $\var(Z)$ is strict positive here):

\begin{center}
\vskip6pt
\begin{tabular}{c|c|c|c} \hline
State $(Z)$ & $z_1$         & $\cdots$ & $z_{m_0}$ \\ \hline
Prob.  & $\frac{S_1}{S}$   & $\cdots$ & $\frac{S_{m_0}}{S}$  \\ \hline
\end{tabular}
\vskip6pt
\end{center}

Thus, the matrix $F_0$ is negative definite, $\utwi{1}\utwi{1}^T+I$ is positive definite,
by Duan and Patton (1998), $F$ is stable.
 Applying Taylor expansion
 to $h(\theta)$ at a point $\theta_*$, we have 
 \[
 \|h(\theta)-F(\theta-\theta_*)\| \leq  c \|\theta-\theta_*\|^2,
 \]
 for some value $c>0$.  Therefore, $(A_2)$ is satisfied by both SAMC and Pop-SAMC.